\documentclass[11pt]{article}

\usepackage{graphicx}
\usepackage{latexsym,amsmath,amsfonts,amscd, amsthm, dsfont}
\usepackage{bm,color}
\usepackage{epsfig,verbatim,epstopdf,graphics}
\usepackage{subfigure}
\usepackage{changebar}
\usepackage{multirow}
\usepackage{morefloats}

\usepackage{algorithmic}

\usepackage{yhmath}
 \usepackage{booktabs} 
 \usepackage{tikz}
\usepackage{verbatim}
\usetikzlibrary{arrows,backgrounds,snakes,shapes}
 \numberwithin{equation}{section}

\graphicspath{{./}{./figure/}}
\allowdisplaybreaks

\newtheoremstyle{plainNoItalics}{}{}{\normalfont}{}{\bfseries}{.}{ }{}

\theoremstyle{plain}
\newtheorem{thm}{Theorem}[section]

\theoremstyle{plainNoItalics}

\newtheorem{exa}[thm]{Example}

\newcommand{\bE}{{\bf E}}

\newcommand{\be}{\begin{eqnarray}}
\newcommand{\ee}{\end{eqnarray}}
\newcommand{\beno}{\begin{eqnarray*}}
\newcommand{\eeno}{\end{eqnarray*}}

\makeatletter

\newcommand{\Rmnum}[1]{\expandafter\@slowromancap\romannumeral #1@}
\makeatother

\begin{document}

\baselineskip=1.8pc


\begin{center}
{\bf
A high order semi-Lagrangian discontinuous Galerkin method for the two-dimensional incompressible Euler equations and the guiding center Vlasov model without operator splitting
}
\end{center}

\vspace{.2in}
\centerline{
Xiaofeng Cai\footnote{
 Department of Mathematical Sciences, University of Delaware, Newark, DE, 19716. E-mail: xfcai@udel.edu.
},
 Wei Guo\footnote{
Department of Mathematics and Statistics, Texas Tech University, Lubbock, TX, 70409. E-mail:
weimath.guo@ttu.edu. Research
is supported by NSF grant NSF-DMS-1620047.
},
Jing-Mei Qiu\footnote{Department of Mathematical Sciences, University of Delaware, Newark, DE, 19716. E-mail: jingqiu@udel.edu. Research of first and last author is supported by NSF grant NSF-DMS-1522777, Air Force Office of Scientific Computing FA9550-12-0318 and University of Delaware.}
}

\bigskip
\noindent
{\bf Abstract.}
In this paper, we generalize a high order semi-Lagrangian (SL) discontinuous Galerkin (DG) method for multi-dimensional linear transport equations without operator splitting developed in Cai et al. (J. Sci. Comput. 73: 514-542, 2017) to the 2D time dependent incompressible Euler equations in the vorticity-stream function formulation and the guiding center Vlasov model. We adopt a local DG method for Poisson's equation of these models. For tracing the characteristics, we adopt a high order characteristics tracing mechanism  based on a prediction-correction technique. The SLDG with large time-stepping size might be subject to extreme distortion of upstream cells. To avoid this problem, we propose a novel adaptive time-stepping strategy by controlling the relative deviation of areas of upstream cells.

\vfill

\newpage

\section{Introduction}

In this paper, we propose a class of high order semi-Lagrangian discontinuous Galerkin (SLDG) methods for the two-dimensional (2D) time dependent incompressible Euler equation in the vorticity stream-function formulation and the guiding center Vlasov model. This is a continuation of our previous research effort on the development of high order non-splitting SLDG methods for 2D linear transport equations \cite{cai2016high} and the Vlasov-Poisson (VP) system \cite{cai2018high}.

The 2D time dependent incompressible Euler equations in the vorticity-stream function formulation reads
\begin{equation}
\begin{split}
\omega_t + \nabla\cdot (\mathbf{u}\omega) = 0, \\
\Delta \Phi = \omega,\
\mathbf{u}
=
( -\Phi_y , \Phi_x),
\end{split}
\label{Euler}
\end{equation}
where $\mathbf{u}$ is the velocity field, $\omega$ is the vorticity of the fluid, and $\psi$ is the stream-function determined by Poisson's equation. The other closely related model concerned in this paper is the guiding center approximation of the 2D Vlasov model, which describes a highly magnetized plasma in the
transverse plane of a tokamak \cite{shoucri1981two,crouseilles2009conservative,frenod2015long,yang2014conservative} and is given as follows,
\begin{align}\label{guiding_center}
\rho_t + \nabla\cdot (\mathbf{E}^\perp\rho) = 0, \\
-\Delta \Phi = \rho,\ \bE^\perp = ( -\Phi_y , \Phi_x)
, \label{poisson}
\end{align}
where $\rho$ is the charge density of the plasma and $\bE$ determined by $\bE = - \nabla \Phi$ is the electric field. We denote $\mathbf{E}=(E_1,E_2)$.
Despite their different application backgrounds, the above two models indeed have an equivalent mathematical formulation up to a sign difference in Poisson's equation. Many research efforts have been devoted to the development of effective numerical schemes for solving the two models. In context of the incompressible model in the vorticity stream-function formulation, we mention the compact finite difference scheme \cite{weinan1996essentially}, the continuous finite element method \cite{liu2001simple}, and the DG method \cite{liu2000high}.
It is worth noting that such a vorticity stream-function formulation is attractive in both theoretical study as well as numerical scheme development for incompressible fluid models. One immediate advantage is that the incompressibility of the velocity field is automatically satisfied without additional divergence cleaning techniques. Meanwhile, this formulation introduces complication of imposing numerical boundary conditions when the viscosity terms are present \cite{bonaventura2018fully,russo1993deterministic,souli1996vorticity,weinan1996vorticity}. We do not pursue this direction and assume periodic boundary conditions in this paper. In the context of the guiding center model, we mention the SL schemes \cite{qiu_shu_sl, crouseilles2009conservative, yang2014conservative}.

In this paper, we propose a high order, stable and efficient numerical scheme for \eqref{Euler} and \eqref{guiding_center} under the DG framework. DG framework is well-known not only for its high order accuracy and ability to resolve fine scale structures, but also for its excellent conservation property, superior performance in long time wave-like simulations, and convenience for hp-adaptive and parallel implementation \cite{cockburn2001runge}.
However, it is well-known that the DG scheme coupled with an explicit Runge-Kutta (RK) time integrator suffers from a stringent CFL time step restriction for stability, despite its many appealing properties such as simplicity for implementation \cite{cockburn2001runge,cai2016high}. Such a drawback becomes more pronounced when the RKDG scheme is applied to \eqref{Euler}. More specifically,
as mentioned in \cite{liu2000high}, the computational cost of the scheme is largely dominated by the Poisson solver, also see the performance study in Section 3. For the RKDG scheme, excessively small time steps have to be chosen for stability; consequently a large number of the Poisson solver will be called in time evolution, leading to immense computational cost. On the other hand, the SL approach is known to be free of the CFL time step restriction by building in the characteristics tracing mechanism in scheme formulation. In this paper, we leverage SL approach to alleviate the efficiency issue associated with the RKDG scheme.

In \cite{cai2016high,cai2018high}, we formulated a class of high order conservative SLDG schemes for solving 2D transport problems with application to the VP system. To the authors' best knowledge, such a method is the first SLDG scheme in the literature that is high order accurate (up to third order accurate), unconditionally stable, mass conservative and free of splitting error for 2D transport simulations. In this work, we consider generalizing the SLDG scheme to solving \eqref{Euler} and \eqref{guiding_center}. The efficiency of such a scheme is realized by taking large time step evolution without any stability issue, while the accuracy is not much compromised. This is very desired when solving \eqref{Euler} and \eqref{guiding_center}, since a much smaller number of calls of the Poisson solver are needed compared with the RKDG scheme, resulting in great computational savings. To  accurately trace the characteristics in a non-splitting fashion, we propose to incorporate a high order two-stage multi-derivative predictor-corrector algorithm proposed in \cite{xiong}. We would like to remark that many existing SL methods for solving  \eqref{Euler} and \eqref{guiding_center} are based on the dimensional splitting approach \cite{qiu_shu_sl, crouseilles2009conservative}. However, unlike the VP system, in the splitting setting it is not straightforward to enhance the splitting error accuracy beyond first order, since the characteristics of the system \eqref{Euler} or \eqref{guiding_center} are more sophisticated and thus more complicated to trace accurately when the transport equation is split. In our earlier work \cite{christlieb2014high}, the integral deferred correction approach is employed to correct splitting errors for a class of high order splitting SL schemes. However, time step constraint due to numerical stability is introduced which impedes efficiency of the SL approach. A detailed comparison on the performance of splitting and non-splitting SL schemes for solving \eqref{Euler} and \eqref{guiding_center} will be conducted in our forthcoming paper. There exist several non-splitting SL schemes in the literature, see \cite{xiong,yang2014conservative}, but they cannot conserve the total mass of the system. Another key ingredient of the proposed scheme is a novel adaptive time-stepping algorithm. By carefully tracking scheme's ability in preserving areas of upstream cells, we are able to adaptively adjust time step sizes to ensure uniformly good approximations to shapes of upstream cells. Numerical evidences in Section 3 show that this adaptive algorithm is very effective in enhancing robustness of the SLDG scheme and removing spurious oscillation induced by unphysical distortion of upstream cells.

The rest of this paper is organized as follows. In Section 2, we formulate the SLDG scheme for solving the guiding center Vlasov model. In particular, three main ingredients consisting of the SLDG framework, a high order characteristics tracing algorithm, and an adaptive time-stepping strategy are introduced. In Section 3, a collection of numerical examples are presented, and schemes' performance under different configurations are evaluated. In Section 4, we conclude the paper with some remarks on future work.

%
%
%
%

\section{Multi-dimensional SLDG algorithm for the nonlinear guiding center Vlasov model}

In this section, we describe our proposed scheme for the 2D guiding center Vlasov model problem. Note that a similar algorithm can be formulated for the 2D incompressible Euler model in vorticity stream-function formulation as well. We start by reviewing the high order truly multi-dimensional SLDG framework originally proposed in \cite{cai2016high} in a linear setting. Then we describe how to incorporate the high order characteristics tracing scheme proposed in \cite{xiong}  in the same SLDG framework for the nonlinear model problem.
Lastly, we propose an adaptive time-stepping strategy, using relative deviation of areas of upstream cells as an adaptive indicator, that greatly improve robustness of the SLDG algorithm in a nonlinear setting.

\subsection{SLDG algorithm framework}
\begin{figure}[h]
\centering
\subfigure[]{
\begin{tikzpicture}
    \draw[black,thin] (0,0.5) node[left] {} -- (5.5,0.5)
                                        node[right]{};
    \draw[black,thin] (0,2.) node[left] {$$} -- (5.5,2)
                                        node[right]{};
    \draw[black,thin] (0,3.5) node[left] {$$} -- (5.5,3.5)
                                        node[right]{};
    \draw[black,thin] (0,5 ) node[left] {$$} -- (5.5,5)
                                        node[right]{};
    \draw[black,thin] (0.5,0) node[left] {} -- (0.5,5.5)
                                        node[right]{};
    \draw[black,thin] (2,0) node[left] {$$} -- (2,5.5)
                                        node[right]{};
    \draw[black,thin] (3.5,0) node[left] {$$} -- (3.5,5.5)
                                        node[right]{};
    \draw[black,thin] (5,0) node[left] {$$} -- (5,5.5)
                                        node[right]{};
    \fill [blue] (3.5,3.5) circle (2pt) node[] {};
    \fill [blue] (5,3.5) circle (2pt) node[] {};
    \fill [blue] (3.5,5) circle (2pt) node[below right] {$A_j$} node[above left] {$c_4$};
    \fill [blue] (5,5) circle (2pt) node[] {};

     \draw[thick,blue] (3.5,3.5) node[left] {} -- (3.5,5)
                                        node[right]{};
      \draw[thick,blue] (3.5,3.5) node[left] {} -- (5,3.5)
                                        node[right]{};
       \draw[thick,blue] (3.5,5) node[left] {} -- (5,5)
                                        node[right]{};
        \draw[thick,blue] (5,3.5) node[left] {} -- (5,5)
                                        node[right]{};
    \fill [red] (1.,1) circle (2pt) node[above right,black] {};
    \fill [red] (3,1) circle (2pt) node[] {};
    \fill [red] (1,2.5) circle (2pt) node[below right] {$A_j^\star$} node[above left] {$c_4^\star$};
    \fill [red] (2.5,2.5) circle (2pt) node[] {};

     \draw[-latex,dashed](3.5,5)node[right,scale=1.0]{}
        to[out=240,in=70] (1,2.50) node[] {};

     \draw (0.5+0.01,2-0.01) node[fill=white,below right] {$A_l$};

     \draw [red,thick] (1,1)node[right,scale=1.0]{}
        to[out=20,in=150] (2,0.7) node[] {};

        \draw [red,thick] (2,0.7)node[right,scale=1.0]{}
        to[out=330,in=240] (3,1) node[] {};
             \draw [red,thick] (1,2.5)node[right,scale=1.0]{}
        to[out=310,in=90] (1.1,2) node[] {};
        \draw [red,thick] (1.1,2)node[right,scale=1.0]{}
        to[out=270,in=80] (1,1) node[] {};

        \draw [red,thick] (1,2.5)node[right,scale=1.0]{}
        to[out=10,in=180] (2.5,2.5) node[] {};

        \draw [red,thick] (3,1)node[right,scale=1.0]{}
        to[out=80,in=280] (2.5,2.5) node[] {};
\end{tikzpicture}
}
\subfigure[]{

\begin{tikzpicture}[scale = 1.3]
    \draw[black,thin] (0,0.5) node[left] {} -- (4,0.5)
                                        node[right]{};
    \draw[black,thin] (0,2.) node[left] {$$} -- (4,2)
                                        node[right]{};
    \draw[black,thin] (0,3.5) node[left] {$$} -- (4,3.5)
                                        node[right]{};
    \draw[black,thin] (0.5,0) node[left] {} -- (0.5,4)
                                        node[right]{};
    \draw[black,thin] (2,0) node[left] {$$} -- (2,4)
                                        node[right]{};
    \draw[black,thin] (3.5,0) node[left] {$$} -- (3.5,4)
                                        node[right]{};

    \fill [red] (1.,1) circle (2pt) node[above right,black] {$A_{j,l}^{n,(\tau)}$};
    \fill [red] (3,1) circle (2pt) node[] {};
    \fill [red] (1,2.5) circle (2pt) node[below right] {$A_j^{n,(\tau)}$} node[above left] {};
    \fill [red] (2.5,2.5) circle (2pt) node[] {};

     \draw (0.5+0.01,2-0.01) node[fill=white,below right] {};

\draw [red,thick] (1,1)node[right,scale=1.0]{}
        to[out=20,in=150] (2,0.7) node[] {};

        \draw [red,thick] (2,0.7)node[right,scale=1.0]{}
        to[out=330,in=240] (3,1) node[] {};
             \draw [red,thick] (1,2.5)node[right,scale=1.0]{}
        to[out=310,in=90] (1.1,2) node[] {};
        \draw [red,thick] (1.1,2)node[right,scale=1.0]{}
        to[out=270,in=80] (1,1) node[] {};

        \draw [red,thick] (1,2.5)node[right,scale=1.0]{}
        to[out=10,in=180] (2.5,2.5) node[] {};

        \draw [red,thick] (3,1)node[right,scale=1.0]{}
        to[out=80,in=280] (2.5,2.5) node[] {};
           \draw (0.5+0.01,2-0.01) node[fill=white,below right] {$A_l$};
         \draw[-latex,ultra thick] (1,1)node[right,scale=1.0]{}
        to  (2,1) node[] {};

             \draw[-latex,ultra thick]  (2,1)node[right,scale=1.0]{}
        to (2,2) node[] {};
             \draw[-latex,ultra thick]  (2,2)node[right,scale=1.0]{}
        to (1,2) node[] {};
             \draw[-latex,ultra thick]   (1,2)node[right,scale=1.0]{}
        to (1,1) node[] {};

\draw [thick] (1,1)-- (3,1) node[] {};
 \draw [thick] (1,2.5) -- (1,1) node[] {};
        \draw [thick] (1,2.5)--(2.5,2.5) node[] {};
        \draw [thick] (3,1)--(2.5,2.5) node[] {};
\end{tikzpicture}

}
\caption{Schematic illustration of the SLDG formulation in two dimensions: quadrilateral approximation to a upstream cell.  }
\label{schematic_2d}
\end{figure}
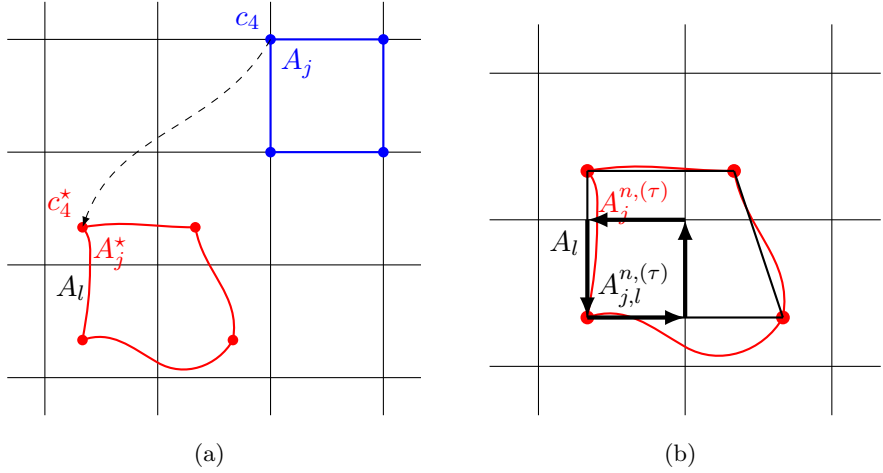

We consider the guiding center Vlasov model \eqref{guiding_center}
on the  2D  domain $\Omega$.
We assume a Cartesian uniform partition of the computational domain $\Omega=\{A_j\}_{j=1}^J$ (see Figure \ref{schematic_2d}) for simplicity.
We define the finite dimensional piecewise polynomial approximation space, $V_h^k = \{ v_h: v_h|_{A_j} \in P^k(A_j) \}$, where $P^k(A_j)$ denotes the set of polynomials of degree at most of $k$ on element $A_j$.
For illustrative purposes, we only present the formulation of the second order SLDG scheme with $P^1$ polynomial space. The generalization, to a third order SLDG scheme with $P^2$ polynomial space and quadratic-curved (QC) quadrilateral approximations to upstream cells, follows a similar procedure discussed in \cite{cai2016high,cai2018high}. The main difference (extra work) involved in a third order SLDG scheme, besides using $P^2$ piecewise polynomials as solution and test function spaces, comes from constructing quadratic curves in approximating sides of upstream cells. Recall that if only regular quadrilaterals are used to approximate upstream cells, then a second order error would be committed, and such an error may become dominant in a nonlinear setting. Numerical evidence will be shown later in the next section in this regard.

In order to update the solution at time level $t^{n+1}$ over the cell $A_j$ based on the solution at time level $t^n$, we employ the weak formulation of characteristic Galerkin method proposed in \cite{Guo2013discontinuous,cai2016high}. Specifically, we consider the following adjoint problem for the time dependent test function $\psi$
\begin{equation}
\psi_t + E_2\psi_x -E_1 \psi_y = 0,\
\text{subject to} \
\psi(t=t^{n+1}) = \Psi(x,y),\
t\in [t^n, t^{n+1} ],
\label{adjoint}
\end{equation}
where $\Psi\in P^k(A_j)$.
The scheme formulation takes advantage of the identity
\begin{equation}
\frac{d}{dt} \int_{  \widetilde{A}_j(t)}  \rho(x,y,t) \psi(x,y,t) dxdy =0,
\end{equation}
where $\widetilde{A}_j(t)$ is a dynamic moving cell, emanating from the Eulerian cell $A_j$ at $t^{n+1}$ backward in time by following characteristics trajectories.
The multi-dimensional  SLDG scheme is formulated as follows:
Given the approximate solution $\rho^n\in V_h^k$ at time $t^n$, find $\rho^{n+1}\in V_h^k$ such that $\forall \Psi\in V_h^k$, we have
\begin{equation}
\int_{A_j} \rho^{n+1} \Psi(x,y) dxdy
=
\int_{A_j^\star} \rho^n \psi(x,y,t^n) dxdy, \quad \mbox{for} \quad j = 1, \cdots, J,
\label{sldg}
\end{equation}
where $\psi$ solves \eqref{adjoint} and $A_j^\star= \tilde{A}_j(t^n)$. $A_j^\star$ is called the upstream cell of $A_j$.
In general, $A_j^\star$ is no longer a rectangle, for example,
see a deformed upstream cell bounded by red curves in Figure \ref{schematic_2d}.
The proposed SLDG method in updating the numerical solution from $\rho^n$ to $\rho^{n+1}$ consists of the following two main steps:
\begin{description}
  \item[1.]  {\bf Construct approximated upstream cells by following characteristics.} 
  Denote four vertices of $A_j$  as $c_q$, with the coordinates $( x_q, y_q )$, $q=1,\cdots,4$.
  We trace characteristics backward in time to $t^n$ for four vertices and then obtain $c_q^\star$ with the new coordinates $( x_q^\star, y_q^\star ),q=1,\cdots,4$. For example, see $c_4$ and $c_4^\star$ in Figure \ref{schematic_2d}.
  Then the upstream cell can be approximated by a quadrilateral determined by four vertices $c_q^\star$.
The new coordinates $( x_q^\star, y_q^\star )$ of $c^\star_q$ are approximated by numerically solving the characteristics equation \eqref{adjoint} in the 2D case, i.e.,
      \begin{equation}
      \begin{cases}
      \frac{dx(t)}{dt} = E_2, \\[3mm]
      \frac{dy(t)}{dt} = -E_1,
       \end{cases}
\quad       \mbox{with} \quad
           \begin{cases}
      x(t^{n+1}) = x_q, \\
      y(t^{n+1} ) = y_q,
     \end{cases}
     \quad
q = 1, 2, 3, 4,
      \label{characteristics}
      \end{equation}
      which is a set of final value problems.
Note that the above equations are non-trivial to solve with high order temporal accuracy, since the $\mathbf{E}$ depends on the unknown $\rho$ via Poisson's equation \eqref{poisson} in a global and nonlinear fashion.
  To circumvent the difficulty, we propose to combine a high order two-stage multi-derivative prediction-correction strategy for tracing characteristics as proposed in \cite{xiong}.
  Such a strategy is described in the context of the proposed SLDG scheme in Section \ref{prediction-correction}.

  \item[2.]  {\bf Update the solution $\rho^{n+1}$ by evaluating the right-hand side of eq.~\eqref{sldg} for $\forall\Psi\in V_h^k$.}
We approximate $A_j^\star$ by a quadrilateral in the previous step. The test function $\psi$ at $t^n$ can be approximated by a polynomial via a least squares procedure by tracking point values of $\psi$ along characteristics.
In order to efficiently evaluate the volume integral in the right-hand side (RHS) of \eqref{sldg}, it is converted into a set of line integrals by the use of Green's theorem. Such an idea is borrowed from CSLAM \cite{lauritzen2010conservative}, and further reformulated in \cite{cai2016high} for the development of an SLDG transport scheme. The above-mentioned procedure is briefly described in Section \ref{section:SLDG}.

\end{description}

\subsection{High order characteristics tracing  algorithm}
\label{prediction-correction}

In this subsection,  we describe a high order predictor-corrector procedure for locating the feet of the characteristics of the guiding center Vlasov model. Such an approach is originally proposed in \cite{qiu2017high}.
It is generalized to the guiding center Vlasov model in \cite{xiong2016conservative}.

We first introduce several shorthand notations. The superscript $^n$ denotes the time level, the superscript $^{(\tau)}$ denotes the formal order of approximation for time discretization, and the subscript $q$ is the index for the vertices of the underlying cell. For example, $(x_q^{n, (\tau)}, y_q^{n, (\tau)})$ is the $\tau$-th order approximation of $(x_q^\star, y_q^\star)$ and $A_j^{n, (\tau)}$ is the quadrilateral determined by the corresponding four vertices.

We will adopt LDG method \cite{arnold2002unified, cockburn1998local, castillo2000priori,zhu2017h} to solve Poisson's equation in each predictor and corrector step. For example, the electric field $\mathbf{E}$ depends on $\rho$ via the Poisson's equation and the time derivative of  $\mathbf{E}$ via another Poisson's equation \eqref{poisson_t} in the corrector step.
Note that the electric field $\mathbf{E}$ is the gradient of potential from the Poisson's equation.
As shown in \cite{arnold2002unified}, using polynomial of degree $k$, the order of convergence for the electric field $\mathbf{E}$ is $k$.
Therefore, if a $(k+1)$-th order of convergence is desired, an LDG scheme with polynomial of degree $k+1$ is needed for solving the Poisson's equation.

The numerical solution $\mathbf{E}_h$ solved by the LDG method are discontinuous across cell boundaries; that is, there are several limits of $\mathbf{E}_h$ from different directions.
For example, $\mathbf{E}_h(x_q^{NE},y_q^{NE}, t^n)$, $\mathbf{E}_h(x_q^{NW},y_q^{NW} , t^n)$, $\mathbf{E}_h(x_q^{SE},y_q^{SE}, t^n)$, $\mathbf{E}_h(x_q^{SW},y_q^{SW} , t^n)$ are not equal, where the superscripts $^{NE}$, $^{NW}$, $^{SE}$, $^{SW}$ are the northeast, northwest, southeast and southwest  limits  of the corresponding functions with respect to $x_q$, respectively. In our implementation, we take the average of $\mathbf{E}_h$ at cell vertices
\begin{equation*}
\mathbf{E}(x_{q},y_{q}, t^n) = \frac{\mathbf{E}_h(x_q^{NE},y_q^{NE}, t^n)+\mathbf{E}_h(x_q^{NW},y_q^{NW} , t^n) +\mathbf{E}_h(x_q^{SE},y_q^{SE}, t^n)+\mathbf{E}_h(x_q^{SW},y_q^{SW} , t^n)}{4}.
\end{equation*}

Next, we present the formulation of a high order predictor-corrector procedure for locating feet of the characteristics in the guiding center Vlasov model.

\begin{description}
\item[First order scheme.]
We start from a first order scheme for tracing characteristics \eqref{characteristics}. We let
\begin{equation}
x_{q }^{n,(1)}
=
x_{q} - E_2( x_q,y_q,t^n ) \Delta t, \
y_{q}^{n,(1)}
=
y_{q} + E_1( x_{q},y_q,t^n) \Delta t, \
\label{eq:first}
\end{equation}
which leads to a first order approximations to $(x_{q}^\star,v_{q}^\star)$.
The  $\mathbf{E}$ depends on $\rho$ at $t^n$ via the Poisson's equation, which can be numerically solved by the LDG method.
Let $A_j^{n,(1)}$ to be the quadrilateral formed by the four upstream vertices $(x_{q }^{n,(1)}, y_{q }^{n,(1)})$, $q=1, 2, 3, 4.$
Then, by the SLDG formulation (to be described in the next subsection)
 \begin{equation}
   \int_{A_j} \rho^{n+1,(1)} \Psi(x,y) dxdy = \int_{ A_j^{n,(1)} } \rho^n \psi(x,y,t^n) dxdy,
   \label{eq: SLDG_t1}
  \end{equation}
we obtain $\rho^{n+1, (1)}$ as a first order approximation in time to $\rho$ at $t^{n+1}$.
Based on $\rho^{n+1, (1)}$, we apply the LDG method to the Poisson's equation \eqref{poisson} again and compute $\mathbf{E}_q^{n+1,(1)}$, which approximates $\mathbf{E}(x_q,y_q, t^{n+1})$ with first order temporal accuracy.

\item[Second order scheme.]
A second order scheme can be built upon the first order one. First,  let
\begin{equation}
\begin{split}
x_{q }^{n,(2)}
=
x_{q} - \frac12\left( E_{2,q}^{n+1,(1)}  + E_2( x_q^{n,(1)} , y_q^{n,(1)} , t^n)  \right)\Delta t, \\
y_{q}^{n,(2)}
=
y_{q} + \frac12 \left(  E_{1,q}^{n+1,(1)} + E_1( x_q^{n,(1)} , y_q^{n,(1)} , t^n )  \right) \Delta t,
\end{split}
\label{eq:second}
\end{equation}
which gives a second order approximation to  $(x_q^\star,y_q^\star)$.
Then the second order approximation solution $\rho^{n+1,(2)}$ is obtained from the SLDG formulation
 \begin{equation}
   \int_{A_j} \rho^{n+1,(2)} \Psi(x,y) dxdy = \int_{ A_j^{n,(2)} } \rho^n \psi(x,y,t^n) dxdy.
   \label{eq: SLDG_t2}
  \end{equation}
Based on $\rho^{n+1,(2)}$, we are able to compute $\mathbf{E}_q^{n+1,(2)}$ from Poisson's equation, which approximates $\mathbf{E}(x_q,y_q , t^{n+1})$ with second order temporal accuracy.

\item[Third order scheme.]
A third order scheme can be designed based on the above second order approximation. Let
\begin{align}
x_{q }^{n,(3)}
=
x_{q} -
E_{2,q}^{n+1,(2)} \Delta t
+\frac{\Delta t^2}{2}
\left(
\frac23( \frac{d}{dt}E_2(x_{q} ,y_q , t^{n+1} ) )^{(2)} + \frac13 \frac{d}{dt} E_2( x_{q}^{n,(2)}, y_{q}^{n,(2)}, t^n )
  \right),
\label{eq:third1}\\
y_{q}^{n,(3)}
=
y_{q}
+
E^{n+1,(2)}_{1,q} \Delta t
-
\frac{\Delta t^2}{2}
\left(
\frac23( \frac{d}{dt}E_1(x_{q} , y_q , t^{n+1} ) )^{(2)} + \frac13 \frac{d}{dt} E_1( x_{q}^{n,(2)},y_{q}^{n,(2)}, t^n )
\right),
\label{eq:third2}
\end{align}
where  $\frac{d}{dt}$ is the material derivative along the characteristic curve, i.e.,
\begin{equation}
\frac{d}{dt}E_s = \frac{\partial E_s}{\partial t} + \frac{\partial E_s}{\partial x} E_2
- \frac{ \partial E_s }{ \partial y } E_1,\ s=1,2.
\label{material}
\end{equation}
Note that on the RHS of the equation \eqref{material}, the partial derivatives are not explicitly given.
The spatial derivative terms $\frac{\partial E_s}{\partial x}, \frac{\partial E_s}{\partial y}$, $s=1,2$ can be approximated by high order DG spatial approximations, while the time derivative term $\frac{\partial E_s}{\partial t}$ can be approximated by utilizing the Vlasov equation (in a Lax-Wendroff spirit in transforming time derivatives into spatial derivatives).
In particular, taking partial time derivative of the 2D Poisson's equation gives
\begin{equation}
\Delta \Phi_t =  (E_2\rho)_x - (E_1\rho)_y.
\label{poisson_t}
\end{equation}
After obtaining $\mathbf{E}$ by solving the original Poisson's equation \eqref{poisson},
the RHS of \eqref{poisson_t} can be constructed by the DG aproximation.
Then we can solve \eqref{poisson_t} by LDG method to get $\frac{\partial \mathbf{E}}{\partial t} = -( (\Phi_t)_x, (\Phi_t)_y ).$
It can be checked by a local truncation error analysis that $(x_{q }^{n,(3)}, y_{q }^{n,(3)})$ is a third order approximation to $(x_q^\star,y_q^\star)$ \cite{xiong2016conservative}. Consequently, the third order approximation solution
$\rho^{n+1,(3)}$ is updated from the SLDG formulation
\begin{equation}
\int_{A_j} \rho^{n+1,(3)} \Psi(x,y) dxdy = \int_{ A_j^{n,(3)} } \rho^n \psi(x,y,t^n) dxdy.
\label{eq: SLDG_t3}
\end{equation}

\end{description}

%
%

\subsection{A two-dimensional SLDG method with quadrilateral upstream cells.}
\label{section:SLDG}

Below, we present the procedure in evaluating the integral $\int_{A_j^\star} \rho^n \psi(x, y, t^n) dxdy$ with a quadrilateral upstream cell $A_j^\star$. In the algorithm design, we have to pay attention to the following two observations, see \cite{cai2016high}.
\begin{itemize}
\item $\Psi = \psi(x,y, t^{n+1})$ is chosen to be polynomial basis functions on $V_h^k$, while, in general $\psi(x, y, t^n)$ is no longer a polynomial. A polynomial function constructed by a least squares procedure is used to approximate $\psi(x, y, t^n)$.
\item  Over the upstream cell $A_j^\star$ (or its approximation $A_{j}^{n,(\tau)}$), $\rho^n(x, y, t^n)$ is discontinuous across Eulerian cell boundaries, see the background Eulerian grid lines in Figure~\ref{schematic_2d}. To properly evaluate the volume integral, one has to perform the evaluation in a sub-area by sub-area manner. Meanwhile, direct evaluation of volume integrals over these irregular-shape sub-areas is very involved in implementation. The proposed strategy is to convert each volume integral into line integrals using Green's Theorem \cite{cai2016high}.
\end{itemize}
Based on these observations, the proposed algorithm consists of two main components. One is the search algorithm that finds the boundaries for each sub-area, i.e. the overlapping region between the upstream cell and background Eulerian cells. The other is the use of Green's theorem that enables us to convert the volume integral to line integrals based on the result of the search algorithm. Below we describe the detailed procedure in evaluating the volume integral over an approximation of upstream cell $A_{j}^{n,(\tau)}$ for the SLDG scheme with $P^1$ polynomial space. Note that the superscript $(\tau)$ is for the order of temporal approximation in the previous subsection.

\begin{figure}[h]
\centering

\subfigure[]{
\begin{tikzpicture}[scale = 1.1]
    \draw[black,thin] (0,0.5) node[left] {} -- (4,0.5)
                                        node[right]{};
    \draw[black,thin] (0,2.) node[left] {$$} -- (4,2)
                                        node[right]{};
    \draw[black,thin] (0,3.5) node[left] {$$} -- (4,3.5)
                                        node[right]{};
    \draw[black,thin] (0.5,0) node[left] {} -- (0.5,4)
                                        node[right]{};
    \draw[black,thin] (2,0) node[left] {$$} -- (2,4)
                                        node[right]{};
    \draw[black,thin] (3.5,0) node[left] {$$} -- (3.5,4)
                                        node[right]{};

    \fill [red] (1.,1) circle (2pt) node[above right,black] {};
    \fill [red] (3,1) circle (2pt) node[] {};
    \fill [red] (1,2.5) circle (2pt) node[below right] {$A_j^{n,(\tau)}$} node[above left] {};
    \fill [red] (2.5,2.5) circle (2pt) node[] {};
   \usetikzlibrary{shapes.geometric}
  \node[fill,star,star points=4, star point ratio=.2] at (2,1) {};
  \node[fill,star,star points=4, star point ratio=.2] at (2,2.5) {};
  \node[fill,star,star points=4, star point ratio=.2] at (1,2) {};
  \node[fill,star,star points=4, star point ratio=.2] at (2.65,2) {};

     \draw (0.5+0.01,2-0.01) node[fill=white,below right] {};
\draw [red,thick] (1,1)node[right,scale=1.0]{}
        to[out=20,in=150] (2,0.7) node[] {};

        \draw [red,thick] (2,0.7)node[right,scale=1.0]{}
        to[out=330,in=240] (3,1) node[] {};
             \draw [red,thick] (1,2.5)node[right,scale=1.0]{}
        to[out=310,in=90] (1.1,2) node[] {};
        \draw [red,thick] (1.1,2)node[right,scale=1.0]{}
        to[out=270,in=80] (1,1) node[] {};

        \draw [red,thick] (1,2.5)node[right,scale=1.0]{}
        to[out=10,in=180] (2.5,2.5) node[] {};

        \draw [red,thick] (3,1)node[right,scale=1.0]{}
        to[out=80,in=280] (2.5,2.5) node[] {};
\draw [thick] (1,1)-- (3,1) node[] {};
 \draw [thick] (1,2.5) -- (1,1) node[] {};
        \draw [thick] (1,2.5)--(2.5,2.5) node[] {};
        \draw [thick] (3,1)--(2.5,2.5) node[] {};

\end{tikzpicture}

}
\subfigure[]{
\begin{tikzpicture}[scale = 1.1]
    \draw[black,thin] (0,0.5) node[left] {} -- (4,0.5)
                                        node[right]{};
    \draw[black,thin] (0,2.) node[left] {$$} -- (4,2)
                                        node[right]{};
    \draw[black,thin] (0,3.5) node[left] {$$} -- (4,3.5)
                                        node[right]{};
    \draw[black,thin] (0.5,0) node[left] {} -- (0.5,4)
                                        node[right]{};
    \draw[black,thin] (2,0) node[left] {$$} -- (2,4)
                                        node[right]{};
    \draw[black,thin] (3.5,0) node[left] {$$} -- (3.5,4)
                                        node[right]{};

    \fill [red] (1.,1) circle (2pt) node[above right,black] {};
    \fill [red] (3,1) circle (2pt) node[] {};
    \fill [red] (1,2.5) circle (2pt) node[below right] {} node[above left] {};
    \fill [red] (2.5,2.5) circle (2pt) node[] {};

     \draw (0.5+0.01,2-0.01) node[fill=white,below right] {};

\draw [thick] (1,1)-- (3,1) node[] {};
 \draw [thick] (1,2.5) -- (1,1) node[] {};
        \draw [thick] (1,2.5)--(2.5,2.5) node[] {};
        \draw [thick] (3,1)--(2.5,2.5) node[] {};
   \usetikzlibrary{shapes.geometric}
  \node[fill,star,star points=4, star point ratio=.2] at (2,1) {};
  \node[fill,star,star points=4, star point ratio=.2] at (2,2.5) {};
  \node[fill,star,star points=4, star point ratio=.2] at (1,2) {};
  \node[fill,star,star points=4, star point ratio=.2] at (2.65,2) {};
  \draw [-latex] (1-0.2,2) -- node[right=3pt]{$\mathcal{L}_q$}(1-0.2,1) node[] {};
  \draw [-latex] (1-0.2,2.5) -- (1-0.2,2) node[] {};

  \draw [-latex] (2.5,2.5+0.2) -- (2,2.5+0.2);
  \draw [-latex] (2,2.5+0.2) -- (1,2.5+0.2);

  \draw [-latex] (1,1-0.2) -- (2,1-0.2) node[] {};
  \draw [-latex] (2,1-0.2) -- (3,1-0.2) node[] {};
  \draw [-latex] (3+0.2,1) -- (2.65+0.2,2) node[] {};
  \draw [-latex] (2.65+0.2,2) -- (2.5+0.2,2.5) node[] {};



\end{tikzpicture}

}
\subfigure[]{
\begin{tikzpicture}[scale = 1.1]
%

\node [below right,blue] at (2,1) {$s_1$};
\node [above,blue] at (2,2.5) {$s_2$};
\node [below,blue] at (0.9,2.) {$s_3$};
\node [above right,blue] at (2.6,2) {$s_4$};
\node [below right, blue] at (2,2) {$c_1$};
    \draw[black,thin] (0,0.5) node[left] {} -- (4,0.5)
                                        node[right]{};
    \draw[black,thin] (0,2.) node[left] {$$} -- (4,2)
                                        node[right]{};
    \draw[black,thin] (0,3.5) node[left] {$$} -- (4,3.5)
                                        node[right]{};
    \draw[black,thin] (0.5,0) node[left] {} -- (0.5,4)
                                        node[right]{};
    \draw[black,thin] (2,0) node[left] {$$} -- (2,4)
                                        node[right]{};
    \draw[black,thin] (3.5,0) node[left] {$$} -- (3.5,4)
                                        node[right]{};

    \fill [red] (1.,1) circle (2pt) node[above right,black] {};
    \fill [red] (3,1) circle (2pt) node[] {};
    \fill [red] (1,2.5) circle (2pt) node[below right] {} node[above left] {};
    \fill [red] (2.5,2.5) circle (2pt) node[] {};

   \usetikzlibrary{shapes.geometric}
  \node[fill,star,star points=4, star point ratio=.2,blue] at (2,1) {};
  \node[fill,star,star points=4, star point ratio=.2,blue] at (2,2.5) {};
  \node[fill,star,star points=4, star point ratio=.2,blue] at (1,2) {};
  \node[fill,star,star points=4, star point ratio=.2,blue] at (2.65,2) {};
  \node[fill,star,star points=4, star point ratio=.2,blue] at (2,2) {};

  \draw [-latex] (2-0.1,2+0.1) -- (2-0.1,2.5-0.1) node[] {};
  \draw [-latex] (2+0.1,2.5-0.1)--(2+0.1,2+0.1 )  node[] {};

   \draw [-latex] (2-0.1,1+0.1) --node[auto]{$\mathcal{S}_q$} (2-0.1,2-0.1) node[] {};
  \draw [-latex] (2+0.1,2-0.1 )--(2+0.1,1+0.1 )  node[] {};
  \draw [-latex] (2.65-0.1,2-0.1) --(2+0.1,2-0.1)  node[] {};
  \draw [-latex] (2+0.1,2+0.1)--(2.65-0.1,2+0.1)  node[] {};
    \draw [-latex] (1+0.1,2+0.1) --(2-0.1,2+0.1)  node[] {};
  \draw [-latex] (2-0.1,2-0.1)--(1+0.1,2-0.1)  node[] {};

     \draw (0.5+0.01,2-0.01) node[fill=white,below right] {};

\draw [thick] (1,1)-- (3,1) node[] {};
 \draw [thick] (1,2.5) -- (1,1) node[] {};
        \draw [thick] (1,2.5)--(2.5,2.5) node[] {};
        \draw [thick] (3,1)--(2.5,2.5) node[] {};
\end{tikzpicture}

}
\caption{Schematic illustration of the search algorithm. }
\label{schematic_search}
\end{figure}
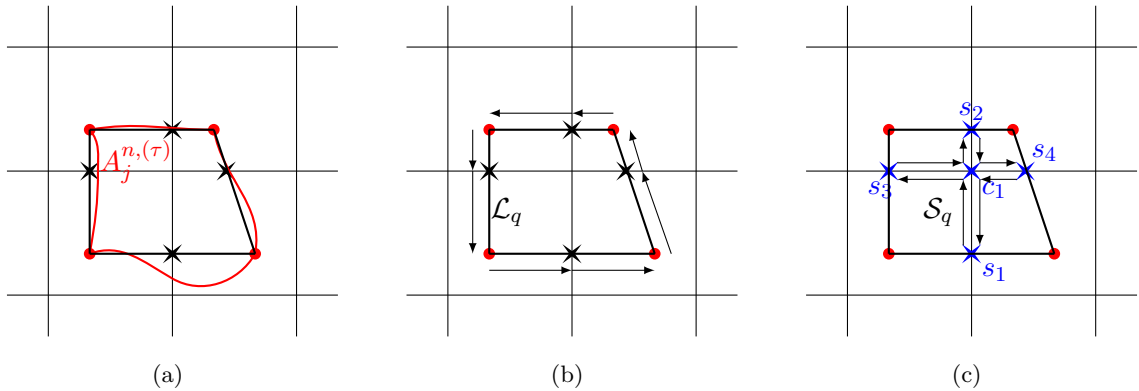

\begin{description}
   \item[(1)] \textbf{\emph{Least squares approximation of test function $\psi(x, y, t^n)$.}}
              Based on the fact that the solution of the adjoint problem \eqref{adjoint} stays unchanged along characteristics, we have
              \begin{equation*}
              \psi( x_q^{n,(\tau)}  ,  y_q^{n,(\tau)}   , t^n ) = \Psi( x_q, y_q ),\ \ q=1,2,\cdots,4.
              \end{equation*}
              Thus, we can reconstruct a unique linear function $\psi^\star(x,y)$ by a least squares strategy that approximates $\psi(x,y,t^n)$ on $A_{j}^{n,(\tau)}$ .

   \item[(2)] 
   \textbf{\emph{Evaluation of the volume integral.}}
   Denote $A_{j,l}^{n,(\tau)}$ as a non-empty overlapping region between the upstream cell $A_j^{n,(\tau)}$ and the background Eulerian cell $A_l$, i.e., $A_{j,l}^{n,(\tau)} = A_{j}^{n,(\tau)} \cap A_l$,
   see Figure \ref{schematic_2d} (b). Then the volume integral, e.g. RHS of eq.~\eqref{eq: SLDG_t1} with $\tau=1$,  becomes
      \begin{equation}
             \int_{A_j^{n,(\tau)}} \rho(x,y,t^{n} )\psi(x,y,t^{n} ) dxdy
           \approx
          \sum_{l\in \varepsilon_j^{n,(\tau)}}^{  } \int_{A_{j,l}^{n,(\tau)} } \rho(x,y,t^n)\psi^\star(x,y)dxdy,
       \label{temp1}
      \end{equation}
      where $\varepsilon_j^{n,(\tau)} =\{ l| A_{j,l}^{n,(\tau)}  \neq \emptyset \}$ and $\psi^\star(x,y)$ is obtained from the previous step.
   Note that the integrands on the RHS of \eqref{temp1} are piecewise   polynomials.
       By introducing two auxiliary function $P(x,y)$ and $Q(x,y)$ such that
       \begin{equation*}
       -\frac{\partial P }{\partial y } + \frac{\partial Q}{\partial x }  =  \rho(x,y,t^n)\psi^\star(x,y),
       \end{equation*}
   the area integral $ \int_{A_{j,l}^{n,(\tau)} } \rho(x,y,t^n)\psi^\star(x,y)dxdy  $ can be converted into line integrals via Green's theorem, i.e.,
   \begin{equation}
      \int_{A_{j,l}^{n,(\tau)} } \rho(x,y,t^n)\psi^\star(x,y)dxdy = \oint_{\partial A_{j,l}^{n,(\tau)}}  Pdx + Qdy,
      \label{Green}
   \end{equation}
   see Figure \ref{schematic_2d} (b). Note that the choices of $P$ and $Q$ are not unique, but the value of the line integrals is independent of the choices. In the implementation, we follow the same procedure in Section 2.1.2 of \cite{lauritzen2010conservative} when choosing $P$ and $Q$.
    In summary, combining \eqref{temp1} and \eqref{Green}, we have the following
\begin{align}
\int_{A_{j}^{n,(\tau)}   } \rho(x,y,t_{n} )\psi(x,y,t_{n} ) dxdy
=&
\sum_{l\in \varepsilon_j^{n,(\tau)}  }^{  }\int_{A_{j,l}^{n,(\tau)} } \rho(x,y,t_{n} )\psi^\star(x,y ) dxdy   \notag \\
=&
\sum_{l\in \varepsilon_j^{n,(\tau)} }^{  } \oint_{\partial A_{j,l}^{n,(\tau)}   }  Pdx + Qdy \notag \\
=&
\sum_{q=1}^{N_o}
\int_{ \mathcal{L}_q } [P  dx +Q  d y ]  + \sum_{q=1}^{N_i}
\int_{ \mathcal{S}_q } [Pdx +Q d y ].
\label{line}
\end{align}
Note that in the above computation, we have organized the liner integrals into two categories: along outer line segments (see Figure \ref{schematic_search} (b)) and along inner line segments (see Figure \ref{schematic_search} (c)). Line segments can be uniquely determined by two end points, which are intersection points of the four sides of the upstream cell with grid lines.
We compute all intersection points and connect them in a counterclockwise orientation to obtain outer line segments, denoted as $\mathcal{L}_q$, $q=1,\cdots,N_o$, see Figure \ref{schematic_search} (b).
The line segments that are aligned with grid lines and enclosed by $A_j^{n,(\tau)}$ are defined as inner line segments, see Figure \ref{schematic_search} (c). Note that there are two orientations along each inner segment, but the corresponding line integrals have to be evaluated in their own sub-area, given that $f^n$ is discontinuous across a inner line segment. For instance, $\overrightarrow{s_1c_1}$ belongs to the left background cell and $\overrightarrow{c_1s_1}$ belongs to the right background cell.

 \end{description}

\subsection{The adaptive time-stepping algorithm}

The proposed SLDG method can be proved to be stable and accurate under large time-stepping size for a linear transport problem with constant coefficients \cite{qiu2011positivity}. However, in a nonlinear setting, in particular with a very large time-stepping size, the upstream cells could be greatly distorted, and the corresponding quadrilaterals or even quadratic-curved quadrilaterals may not be adequate to approximate actual shapes of upstream cells. Consequently, numerical oscillations will be generated, e.g. see Figure~\ref{figure:vortex_WL} for the vortex patch problem in the next section.

Due to the divergence-free constraint on the electric field of the guiding center Vlasov model, the areas of upstream cells should be preserved, i.e., $\text{area}(A_j)=\text{area}(A_j^\star)$ in Figure \ref{schematic_2d}.
If, at the discrete level, the areas of upstream cells are preserved, the local maximum principle in terms of cell averages will be maintained; if the upstream cells are too distorted, e.g. quadrilateral shapes cannot offer adequate approximations (the area of a numerical upstream cell greatly deviates from the actual area), unphysical numerical oscillations may appear. In this section, we propose to measure the $L^\infty$ norm of relative deviation of the area for upstream cells  and use it as an indicator to adaptively select appropriate time-stepping sizes, thus to improve robustness of the SLDG scheme in a nonlinear setting. Below, we first summarize the coupling between the SLDG framework with a third order characteristics tracing algorithm in the Main Algorithm, followed by a detailed description of the Adaptive Time-Stepping Algorithm.
%

\begin{center}
\textbf{Main Algorithm}
\end{center}
\begin{description}
\item[0.] {\bf Initially,} choose a parameter $CFL_{max}$ and let $irefine=0$.
\item[1.] {\bf The first order prediction:}
   \begin{description}
     \item[1.1] Solve the electric field $\mathbf{E}$ by the LDG method, based on the solution $\rho^n$.
     \item[1.2] Trace the characteristics \eqref{characteristics} for a time step $\Delta t$ by the first order scheme \eqref{eq:first} and construct approximate upstream cells $A_j^{n,(1)}$. Perform the ``Adaptive Time-Stepping Algorithm".
     \item[1.3] Evolve the solution $\rho^n$ by using  SLDG  to get $\rho^{n+1,(1)}$.
   \end{description}
\item[2.] {\bf The second order prediction:}
   \begin{description}
     \item[2.1] Solve the electric field $\mathbf{E}$ by the LDG method, based on the solution $\rho^{n+1,(1)}$.
     \item[2.2] Trace the characteristics \eqref{characteristics} for a time step $\Delta t$ by the second order scheme \eqref{eq:second} and construct approximate upstream cells $A_j^{n,(2)}$. Perform the ``Adaptive Time-Stepping Algorithm".
     \item[2.3] Evolve the solution $\rho^{n+1,(1)}$ by using SLDG to get $\rho^{n+1,(2)}$.
   \end{description}
\item[3.] {\bf The third order correction:}
   \begin{description}
     \item[3.1] Solve the electric field $\mathbf{E}$ by the LDG method, based on the solution $\rho^{n+1,(2)}$.
     \item[3.2] Trace the characteristics \eqref{characteristics} for a time step $\Delta t$ by the second order scheme \eqref{eq:third1}-\eqref{eq:third2} and construct approximate upstream cells $A_j^{n,(3)}$. Perform the ``Adaptive Time-Stepping Algorithm".
     \item[3.3] Evolve the solution $\rho^{n+1,(2)}$ by using SLDG  to get $\rho^{n+1}$.
   \end{description}
\end{description}
%
\begin{center}
\textbf{Adaptive Time-Stepping Algorithm}
\end{center}
\fbox{
\begin{minipage}[htb]{0.9\linewidth}
\begin{itemize}
  \item Compute $\theta= \max_j \frac{ \text{area}\left(A_j^{n,(\tau)}\right) - \text{area}\left(A_j\right)  }{ \text{area}\left(A_j\right)}$.

 Let $\delta_M$ and $\delta_m$ be prescribed thresholds for decreasing and increasing CFL number. In our simulations,  $\delta_M=1\%$ and $\delta_m=0.3\%$.
       \begin{description}
         \item[if] $\theta>\delta_M$, {\bf then}
            we let $CFL=\frac23 CFL$, $irefine = 1$ and go back to Step 1.2.
         \item[else if] $\theta<\delta_m$, $irefine =0$, and $CFL\neq CFL_{\max}$, {\bf then}
         $CFL=\min\{ \frac32 CFL, CFL_{\max} \}$ go back to Step 1.2.
         \item[else] Continue to the next step.
         \item[end if]
       \end{description}
\end{itemize}

\end{minipage}
}

\section{Numerical Results}
\label{numerical}

In this section, for the 2D incompressible Euler equation in vorticity stream-function formulation \eqref{Euler} and the guiding center Vlasov model \eqref{guiding_center}, we examine the performance of the proposed SLDG method with second/third order temporal accuracy, denoted by  SLDG+time2/3, with quadrilateral or quadratic-curved (QC) quadrilateral approximation to upstream cells (using the notation without or with -QC), with $P^k$ local discontinuous Galerkin method (using the notation +$P^k$ LDG), without or with the WENO limiter \cite{zhong2013simple} (using the notation without or with +WL).
In all our numerical tests, we let the time step size
\begin{equation*}
\Delta t = \frac{ CFL }{ a/\Delta x + b/\Delta y },
\end{equation*}
in which $CFL$ is specified for different runs.
For the incompressible Euler equation,
 $a=\max (|u|), b=\max(|v|)$.
For the guiding center Vlasov model,
 $a=\max (|E_2|), b=\max(|E_1|)$.
 For example, $P^2$ SLDG-QC+$P^3$ LDG+time3+WL-CFL3 refers to the SLDG scheme with $P^2$ polynomial space, with quadratic-curved quadrilateral approximation to upstream cells, with $P^3$ LDG scheme in solving Poisson's equation, using the third order scheme in characteristics tracing, with the WENO limiter and $CFL=3$. We also apply the proposed SLDG method with the adaptive time-stepping strategy to improve the robustness and efficiency of the method.

For both models, besides mass conservation, the following physical quantities remain constant over time
  \begin{description}
  \item [1.] Mass:
              \begin{equation*}
 \int_{\Omega} \omega dxdy,  \quad \mbox{(Euler),}
\quad           \int_{\Omega} \rho dxdy, \quad \mbox{(Vlasov).}
            \end{equation*}
   \item[2.] Energy:
            \begin{equation*}
            \| \mathbf{u} \|_{L^2}^2 = \int_{\Omega} \mathbf{u}\cdot\mathbf{u}dxdy, \quad \mbox{(Euler),}
\quad
            \| \mathbf{E} \|_{L^2}^2 = \int_{\Omega} \mathbf{E}\cdot\mathbf{E}dxdy, \quad \mbox{(Vlasov).}
            \end{equation*}
   \item[3.]  Enstrophy:
            \begin{equation*}
            \| \omega  \|_{L^2}^2 = \int_{\Omega} \omega^2 dxdy,  \quad \mbox{(Euler),}
\quad            \| \rho  \|_{L^2}^2 = \int_{\Omega} \rho^2 dxdy, \quad \mbox{(Vlasov).}
            \end{equation*}
 \end{description}
Tracking relative deviations of these quantities numerically provides a good measurement of the quality of numerical schemes.
For our numerical tests shown below, all SLDG schemes can conserve total mass up to the round-off error: $O(10^{-13})$ as expected; while we keep track of energy and enstrophy over time to compare performances of SLDG schemes in various settings. Furthermore, due to the incompressibility constraint of ${\bf u}$ (Euler) or ${\bf E}$ (Vlasov), the area of an upstream cell should be preserved. We also track relative deviations of areas of upstream cells, to better understand how we approximate shapes of upstream cells. In our adaptive time-stepping strategy, we use the relative deviation of areas of upstream cells  as a metric to determine if time-stepping size should be increased, reduced or kept the same. In our simulations, we use the threshold of $0.3\%$ for increasing time-stepping sizes; and the threshold of $1\%$ for reducing time-stepping sizes.

Below we present four benchmark examples to assess and compare performances of SLDG schemes with various configurations. Comparisons are made in terms of numerical errors for smooth problems, CPU time, ability to resolve solution structures, robustness, and performance in conserving physical invariants.
Based on all data we collected, $P^2$ SLDG-QC, using $P^3$ or $P^2$ LDG solver for Poisson's equation, coupled with the third order characteristics tracing scheme and the adaptive time-stepping strategy is considered to be an optimal configuration that well balances its performance in effectiveness, efficiency and robustness. The choice of using $P^3$ or $P^2$ LDG solver for Poisson's equation is a trade-off between accuracy (effectiveness in resolving solutions) and CPU cost. Instead of drawing a definite conclusion, we refer to the efficiency comparison presented in Figure~\ref{figure:CPU} for a smooth test, in Figures~\ref{figure:KH_contour} for performance in resolving solution structures,  as well as in Figures~\ref{figure:KH_area}-\ref{figure:KH_norm} in conserving physical invariants. Finally, we would like to remark that the CFL constraint for a Runge-Kutta DG method is known to be $\frac{1}{2k+1}$ where $k$ is the degree of polynomial space. That is $CFL \le 1/3$ for $P^1$ and $CFL\le1/5$ for $P^2$. By using the SLDG algorithm, while maintaining good resolution of solution structures and preservation of invariants, we are able to take CFL as large as $3$ ($9$ times as large for $P^1$ and $15$ times as large for $P^2$), leading to huge computational savings. Note that the dominant CPU cost per time step is the LDG solver for the Poisson equation (as shown in Table 3.3), the extra CPU cost from SLDG method in characteristics tracing and in evaluation of line integrals, compared with that from a Runge-Kutta DG method, will not play a significant role.



\begin{exa}
\label{exa:1}
(Accuracy and convergence test).
Consider the incompressible Euler equation \eqref{Euler}
on the domain $[0,2\pi]\times[0,2\pi]$ with the initial condition
\begin{equation}
\omega (x,y,0) = -2 \sin(x) \sin(y)
\end{equation}
and periodic boundary conditions.
The exact solution stays stationary as $\omega(x,y,t) = -2 \sin(x) \sin(y)$.
We test the spatial convergence, temporal convergence and CPU cost of the proposed SLDG methods for solving \eqref{Euler} up to time $T=1$.

First,  we test the spatial convergence of the proposed SLDG schemes and summarize results in Table~\ref{table:p1}, \ref{table:p2CFL1}, and \ref{table:p2qc} for $P^1$ SLDG scheme, $P^2$ SLDG scheme and $P^2$ SLDG-QC scheme respectively.  The schemes are coupled with LDG schemes of different orders for solving Poisson's equation and  characteristics tracing schemes of different orders. We let $CFL=1$, for which the spatial error still dominates. Expected orders of convergence are observed for all these different settings. The data reported in Tables~\ref{table:p1}, \ref{table:p2CFL1}, and \ref{table:p2qc} are organized into a CPU time versus error log-log plot in Figure~\ref{figure:CPU} to benchmark performances of SLDG schemes with various configurations. We demonstrate the temporal order of convergence in Table~\ref{table:temporal} by varying $CFL$ numbers. Based on all data collected, we make the following observations.
\begin{enumerate}
  \item \emph{For a $P^{k}$ $k=1,2$ SLDG(-QC) scheme, in order to attain $k$th order accuracy, an LDG scheme with $P^{k+1}$ solution space for Poisson's equation is needed.} In Table \ref{table:p1}, we observe that $P^1$ SLDG with $P^1$ LDG is only first order accurate in $L^\infty$ error and $P^1$ SLDG with $P^2$ LDG is second order accurate in $L^\infty$ error. Both $L^1$ and $L^\infty$ errors become smaller when a $P^2$ LDG scheme is used. Note that the velocity (in Euler) or the electric field (in Vlasov) is the gradient of the potential function solved from Poisson's equation; by taking one order of spatial derivative, the order of convergence becomes one order less \cite{arnold2002unified}. Similarly, in Table \ref{table:p2qc} for $P^2$ SLDG-QC scheme, we observe that $P^2$ SLDG-QC with $P^2$ LDG displays a second order spatial convergence, while $P^2$ SLDG-QC with $P^3$ LDG is third order accurate.
\item {\em For a $P^2$ SLDG scheme (without QC), the spatial convergence is of second order due to the error in approximating upstream cells.} We test the schemes with the $P^2$ and $P^3$ LDG schemes, and with the second and third order characteristics tracing schemes (time2 and time3, respectively) and report results in Table \ref{table:p2CFL1}. We also provide the CPU cost and errors associated with each configuration in the table, while plotting CPU time versus error in Figure~\ref{figure:CPU}. For this example, it seems that a $P^2$ LDG is more cost effective, when coupled with the $P^2$ SLDG scheme for the transport equation.
\item {\em LDG Poisson solver dominates the CPU cost of the simulation.} CPU time for simulations with various configurations is collected. In particular, the ratio of CPU(LDG)/CPU(SLDG) is reported in Table~\ref{table:p2qc} for $P^2$ SLDG-QC schemes. It is observed that the LDG Poisson solver dominates the CPU cost. Base on such observation, we suggest to use $P^2$ SLDG-QC (rather than the $P^2$ SLDG scheme) for better computational performance. Data points reported in Figure~\ref{figure:CPU}  support the same conclusion that $P^2$ SLDG-QC is more cost effective than $P^2$ SLDG. We report the CPU time versus error study in Figure~\ref{figure:CPU}. It is observed that the configuration of $P^2$ SLDG-QC with $P^3$ LDG and third order characteristics tracing scheme is the most efficient, once the error tolerance is below certain threshold.
\item {\em Temporal convergence is observed for CFL ranging from $1$ to $6$.} Table \ref{table:temporal} summarizes the errors and the corresponding temporal convergence rates for $P^2$ SLDG-QC+$P^3$ LDG with first, second and third order characteristic tracing schemes and with $CFL$ ranging from $1$ to $6$. To make the temporal error dominant, we use a spatial mesh of $100\times100$ elements. Expected orders of convergence are observed. Higher order characteristics tracing schemes offer not only better convergence rates (only slightly better rate when comparing second and third order schemes), but also smaller errors. Notice that the third order characteristics tracing scheme would cost about $2.5$ times as much CPU time as the second order one, if other settings are the same. Note that, for the second and third order characteristics tracing schemes, the LDG solver (the subroutine with dominant CPU cost) will be called two and five times, respectively. For example, compare the CPU cost of $P^2$ SLDG+ $P^2$ LDG+time2 and  $P^2$ SLDG+ $P^2$ LDG+time3; and the CPU cost of  $P^2$ SLDG+ $P^3$ LDG+time2 and  $P^2$ SLDG+ $P^3$ LDG+time3 in Table~\ref{table:p2CFL1}. 
\end{enumerate}


\begin{table}[!ht]
\caption{
Example~\ref{exa:1} the incompressible Euler equations.
Errors, orders and CPU times (sec) of $P^1$ SLDG+$P^k$ LDG+time2, $k=1,2$. $T=1$. $CFL=1$.
 }
\centering
\begin{tabular}{ccccc cc c}
\hline
{  Mesh }  &{$L^1$ error} & Order    &{$L^2$ error} & Order  & {$L^\infty$ error} & Order & CPU (sec)\\
\hline

  \multicolumn{8}{l}{ $P^1$ SLDG+$P^1$ LDG+time2  }
\\
    $20^2$ &     1.62E-02 & &     2.34E-02 & &     1.32E-01 & & 0.06\\
    $40^2$ &     4.35E-03 &     1.90 &     7.21E-03 &     1.70 &     6.75E-02 &     0.96 & 0.82\\
    $60^2$ &     2.00E-03 &     1.91 &     3.55E-03 &     1.74 &     5.00E-02 &     0.74 & 3.68 \\
    $80^2$ &     1.15E-03 &     1.93 &     2.13E-03 &     1.79 &     3.81E-02 &     0.95 & 13.85 \\
   $100^2$ &     7.41E-04 &     1.96 &     1.41E-03 &     1.83 &     3.08E-02 &     0.96 & 33.87 \\
    \hline

   \multicolumn{8}{l}{ $P^1$ SLDG+$P^2$ LDG+time2 }
\\

    $20^2$ &     1.17E-02 & &     1.57E-02 & &     8.55E-02 & & 0.21\\
    $40^2$ &     2.94E-03 &     1.99 &     4.00E-03 &     1.97 &     2.48E-02 &     1.78 & 2.47\\
    $60^2$ &     1.31E-03 &     1.99 &     1.79E-03 &     1.98 &     1.16E-02 &     1.87 & 12.68\\
    $80^2$ &     7.45E-04 &     1.97 &     1.02E-03 &     1.96 &     6.71E-03 &     1.91 & 49.12\\
   $100^2$ &     4.75E-04 &     2.02 &     6.49E-04 &     2.01 &     4.34E-03 &     1.95  & 131.11\\

   \hline

\end{tabular}
\label{table:p1}
\end{table}

\begin{table}[!ht]
\caption{
Example~\ref{exa:1} the incompressible Euler equations.
Errors, orders and CPU times (sec) of $P^2$ SLDG+$P^{k}$ LDG with different order temporal accuracy, $k=2,3$. $T=1$. $CFL=1$.
 }
\centering
\begin{tabular}{ccccc cc c}
\hline
{  Mesh }  &{$L^1$ error} & Order    &{$L^2$ error} & Order  & {$L^\infty$ error} & Order & CPU \\
\hline

  \multicolumn{8}{l}{ $P^2$ SLDG+$P^2$ LDG+time2  }
\\
    $20^2$ &     9.36E-03 & &     1.34E-02 & &     1.06E-01 & & 0.22\\
    $40^2$ &     2.23E-03 &     2.07 &     3.30E-03 &     2.03 &     3.06E-02 &     1.80 & 2.46\\
    $60^2$ &     9.49E-04 &     2.10 &     1.42E-03 &     2.08 &     1.39E-02 &     1.94 & 12.66\\
    $80^2$ &     5.63E-04 &     1.81 &     8.46E-04 &     1.79 &     8.31E-03 &     1.79 &46.28\\
   $100^2$ &     3.52E-04 &     2.11 &     5.30E-04 &     2.09 &     5.38E-03 &     1.95 &127.00\\
    \hline

   \multicolumn{8}{l}{ $P^2$ SLDG+$P^3$ LDG+time2 }
\\

    $20^2$ &     6.36E-03 & &     9.00E-03 & &     5.58E-02 & & 0.55\\
    $40^2$ &     1.33E-03 &     2.26 &     1.94E-03 &     2.22 &     1.43E-02 &     1.96 &6.68\\
    $60^2$ &     5.74E-04 &     2.07 &     8.19E-04 &     2.12 &     6.46E-03 &     1.97 &47.53\\
    $80^2$ &     3.23E-04 &     2.00 &     4.74E-04 &     1.90 &     3.72E-03 &     1.92 &147.69\\
   $100^2$ &     1.96E-04 &     2.24 &     2.86E-04 &     2.27 &     2.40E-03 &     1.97 &379.80\\

   \hline

  \multicolumn{8}{l}{ $P^2$ SLDG+$P^2$ LDG+time3  }
\\
    $20^2$ &     5.94E-03 & &     8.55E-03 & &     6.53E-02 & & 0.65\\
    $40^2$ &     1.24E-03 &     2.26 &     1.82E-03 &     2.23 &     1.50E-02 &     2.12 &7.34\\
    $60^2$ &     5.30E-04 &     2.10 &     7.62E-04 &     2.15 &     6.31E-03 &     2.13 &37.41\\
    $80^2$ &     3.03E-04 &     1.95 &     4.40E-04 &     1.91 &     3.49E-03 &     2.06 &130.08\\
   $100^2$ &     1.82E-04 &     2.28 &     2.65E-04 &     2.28 &     2.21E-03 &     2.04 &343.46\\
    \hline

   \multicolumn{8}{l}{ $P^2$ SLDG+$P^3$ LDG+time3 }
\\

    $20^2$ &     5.94E-03 & &     8.35E-03 & &     4.77E-02 & & 1.59\\
    $40^2$ &     1.29E-03 &     2.20 &     1.88E-03 &     2.15 &     1.02E-02 &     2.23 &20.25\\
    $60^2$ &     5.65E-04 &     2.04 &     8.06E-04 &     2.09 &     4.22E-03 &     2.16 &131.56\\
    $80^2$ &     3.19E-04 &     1.99 &     4.69E-04 &     1.88 &     2.46E-03 &     1.87 &427.40\\
   $100^2$ &     1.93E-04 &     2.24 &     2.83E-04 &     2.26 &     1.53E-03 &     2.14 &1000.31\\

   \hline
\end{tabular}
\label{table:p2CFL1}
\end{table}

\begin{table}[!ht]
\caption{
Example~\ref{exa:1} the incompressible Euler equations.
Errors, orders and CPU times (sec) of $P^2$ SLDG-QC+$P^{k}$ LDG+time3, $k=2,3$. $T=1$. $CFL=1$.
 }
\centering
\begin{tabular}{ccccc cc c c}
\hline
{  Mesh }  &{$L^1$ error} & Order    &{$L^2$ error} & Order  & {$L^\infty$ error} & Order & CPU & $\frac{\text{CPU(LDG)}}{\text{CPU(SLDG)} }$\\
\hline

  \multicolumn{9}{l}{ $P^2$ SLDG-QC+$P^2$ LDG+time3  }
\\
    $20^2$ &     4.10E-03 & &     6.32E-03 & &     7.94E-02 & & 0.68&   0.45/0.10 \\
    $40^2$ &     6.07E-04 &     2.76 &     9.80E-04 &     2.69 &     1.76E-02 &     2.17 &7.43 & 5.96/0.75\\
    $60^2$ &     2.26E-04 &     2.43 &     3.65E-04 &     2.43 &     7.54E-03 &     2.10 &38.82 & 33.72/2.74\\
    $80^2$ &     1.19E-04 &     2.23 &     1.92E-04 &     2.25 &     4.15E-03 &     2.08 &139.64 & 128.30/5.93\\
   $100^2$ &     7.37E-05 &     2.15 &     1.20E-04 &     2.11 &     2.63E-03 &     2.05 &384.87 & 362.40/12.28\\
    \hline

   \multicolumn{9}{l}{ $P^2$ SLDG-QC+$P^3$ LDG+time3 }
\\

    $20^2$ &     2.19E-03 & &     2.82E-03 & &     1.36E-02 & & 1.63 & 1.29/0.11 \\
    $40^2$ &     2.71E-04 &     3.01 &     3.56E-04 &     2.99 &     1.93E-03 &     2.81 & 17.63 & 15.41/0.82\\
    $60^2$ &     8.00E-05 &     3.01 &     1.04E-04 &     3.02 &     5.33E-04 &     3.18 & 113.77 & 106.31/2.86\\
    $80^2$ &     3.37E-05 &     3.01 &     4.44E-05 &     2.97 &     2.28E-04 &     2.96 & 338.95 & 322.37/6.45 \\
   $100^2$ &     1.71E-05 &     3.05 &     2.24E-05 &     3.08 &     1.14E-04 &     3.09 & 1003.01 & 971.15/12.88\\
   \hline

\end{tabular}
\label{table:p2qc}
\end{table}

\begin{figure}[h!]
\centering
\includegraphics[width=110mm]{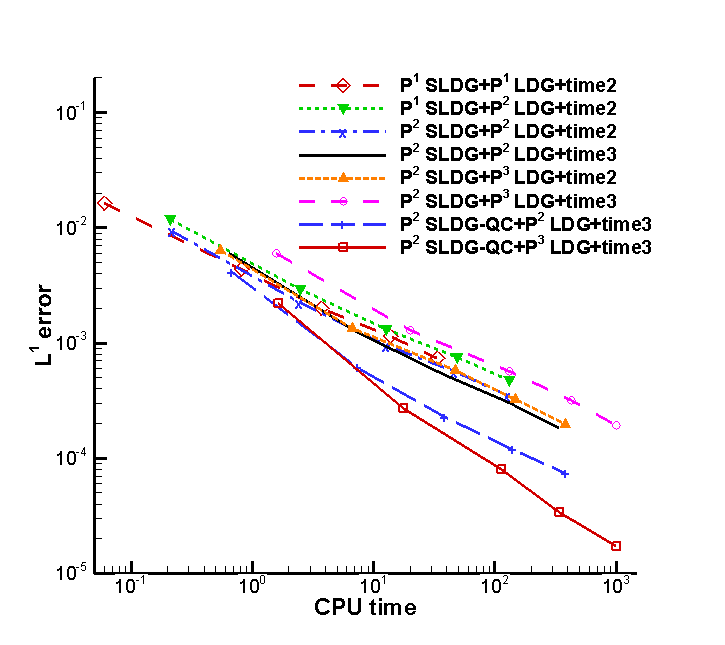}
\caption{Example \ref{exa:1}. $L^{1}$ error versus CPU time (s) in log-log plot. The data points are all from Table~\ref{table:p1}-\ref{table:p2qc}.}
\label{figure:CPU}
\end{figure}

\begin{table}[!ht]
\caption{
The incompressible Euler equations \eqref{Euler} on the domain $[0,2\pi]\times[0,2\pi]$ with the initial condition
$
\omega (x,y,0) = -2 \sin(x) \sin(y).
$
Periodic boundary conditions in two directions.
Temporal order of convergence of $P^2$ SLDG-QC with the mesh of $100\times100$. $T=1$.
The time-stepping size for $P^k$ SLDG-(QC) is $\Delta t =\frac{CFL}{ \frac{\max(|u|) }{\Delta x} + \frac{\max(|v|)}{\Delta y} }$.
 }
\centering
\begin{tabular}{ccccc cc}
\hline
{  $CFL$ }  &{$L^1$ error} & Order    &{$L^2$ error} & Order  & {$L^\infty$ error} & Order \\
\hline

  \multicolumn{7}{l}{ $P^2$ SLDG-QC+$P^3$ LDG+time1  }
\\
   1 &     1.34E-02 & &     1.87E-02 & &     6.55E-02 & \\
   2 &     2.73E-02 &     1.02 &     3.82E-02 &     1.03 &     1.37E-01 &     1.07 \\
   3 &     4.06E-02 &     0.98 &     5.74E-02 &     1.01 &     2.10E-01 &     1.06 \\
   4 &     5.59E-02 &     1.11 &     7.99E-02 &     1.15 &     3.00E-01 &     1.23 \\
   5 &     6.79E-02 &     0.87 &     9.77E-02 &     0.90 &     3.72E-01 &     0.97 \\
   6 &     8.18E-02 &     1.02 &     1.19E-01 &     1.06 &     4.60E-01 &     1.16 \\
    \hline

   \multicolumn{7}{l}{ $P^2$ SLDG-QC+$P^3$ LDG+time2 }
\\

   1 &     3.27E-05 & &     4.27E-05 & &     5.96E-04 & \\
   2 &     1.41E-04 &     2.11 &     1.87E-04 &     2.13 &     9.46E-04 &     0.67 \\
   3 &     3.72E-04 &     2.39 &     5.29E-04 &     2.56 &     2.38E-03 &     2.27 \\
   4 &     7.89E-04 &     2.62 &     1.20E-03 &     2.84 &     5.26E-03 &     2.76 \\
   5 &     1.40E-03 &     2.56 &     2.18E-03 &     2.69 &     9.65E-03 &     2.72 \\
   6 &     2.32E-03 &     2.78 &     3.70E-03 &     2.90 &     1.66E-02 &     2.96 \\

   \hline

   \multicolumn{7}{l}{ $P^2$ SLDG-QC+$P^3$ LDG+time3 }
\\

   1 &     1.71E-05 & &     2.24E-05 & &     1.14E-04 & \\
   2 &     2.63E-05 &     0.62 &     3.35E-05 &     0.58 &     2.04E-04 &     0.84 \\
   3 &     6.06E-05 &     2.06 &     8.45E-05 &     2.28 &     5.58E-04 &     2.48 \\
   4 &     1.21E-04 &     2.41 &     1.70E-04 &     2.42 &     1.01E-03 &     2.06 \\
   5 &     2.05E-04 &     2.36 &     2.78E-04 &     2.22 &     1.11E-03 &     0.44 \\
   6 &     3.50E-04 &     2.93 &     4.72E-04 &     2.91 &     1.59E-03 &     1.95 \\
   \hline
\end{tabular}
\label{table:temporal}
\end{table}

\end{exa}

\begin{exa}
(Kelvin-Helmholtz instability problem).
This example is the 2D guiding center model problem with the initial condition
\begin{equation}
\rho_0(x,y) = \sin(y) + 0.015 \cos(kx)
\end{equation}
and periodic boundary conditions on the domain $[0,4\pi]\times[0,2\pi]$.
We let $k=0.5$, which will create a Kelvin-Helmholtz instability.

We test our schemes in various settings. No WENO limiter is used.  Several representative figures are shown in Figure \ref{figure:KH_contour1}, \ref{figure:KH_contour}, \ref{figure:KH_CFL}, \ref{figure:KH_1d}, \ref{figure:KH_area}, and \ref{figure:KH_norm}.
\begin{enumerate}
  \item \emph{Compare performances of a third order $P^2$ SLDG-QC scheme and a second order $P^1$ SLDG scheme.} In Figure \ref{figure:KH_contour1}, we plot the contour of the solution computed by $P^1$ SLDG with the mesh of $100\times100$ elements as well as the refined mesh of $200\times200$ elements; we also plot the contour of the solution computed by $P^2$ SLDG-QC with the mesh of $100\times100$ elements. By carefully comparing these results, we observe that the third order scheme offers better resolution; and its solution is consistent with that from a second order scheme with the refined mesh ($200\times200$). In Figure \ref{figure:KH_1d}, we compare these solutions via their 1D cuts along the line $y=\pi$. It is observed that the second order solution with the $200\times200$ mesh is more consistent with that of the third order solution with the $100\times100$ mesh.
  \item \emph{$P^2$ SLDG-QC+$P^{3}$ LDG with adaptive $CFL$ versus that with fixed $CFL=3$.} We perform the simulation with adaptive $CFL$, and set the initial $CFL$ to be 3. As the solution evolves, the $CFL$ will be dynamically adjusted according to the adaptive time-stepping algorithm we proposed. In particular, if the $L^\infty$ norm of relative deviation of areas of upstream cells exceeds a threshold ($1\%$), or is below another threshold ($0.3\%$), the time-stepping size will be reduced or increased. The $CFL$ history of the $P^2$ SLDG-QC+$P^3$ LDG+time3 with adaptive $CFL$ is shown in Figure~\ref{figure:KH_CFL}.  Figure \ref{figure:KH_contour} displays the contour plot of the solution with adaptive $CFL$; while Figure \ref{figure:KH_1d} (b) shows the 1D cut of the solutions of $P^2$ SLDG-QC and $P^{3}$ LDG with fixed $CFL=3$ and with adaptive $CFL$. The solution with $CFL=1$ is plotted in the same figure, as a reference solution. The scheme with adaptive $CFL$ is observed to be able to capture the solution well.
  \item {\em SLDG-QC scheme with adaptive CFL: comparison for using $P^2$ or $P^3$ LDG scheme for Poisson's equation.} We find comparable performance of the adaptive schemes using $P^2$ and $P^3$ LDG solving Poisson's equation in resolving solution structures in Figure~\ref{figure:KH_contour} and in preserving upstream cell areas as well as physical invariants in Figure~\ref{figure:KH_area}-\ref{figure:KH_norm}. The scheme with $P^3$ performs only slightly better in preserving physical invariants; however, the CPU cost of the scheme with the $P^3$ LDG is twice as much as that of the same scheme but with $P^2$ LDG, see Table~\ref{table:p2qc}. Note that for the previous smooth example, as shown in Figure~\ref{figure:CPU}, when the error is below certain threshold, $P^3$ LDG scheme is doing slightly better; when the error is above that threshold (in other words, the solution has not been well-resolved), a $P^2$ LDG Poisson solver may better balance efficiency and effectiveness.
 \item {\em Time histories of $L^\infty$ norm of relative deviation of areas of upstream cells are shown in Figure~\ref{figure:KH_area}.} On its subplot (a), we compare the performance of schemes with fixed $CFL=1$, $CFL=3$ and adaptive $CFL$ with initial $CFL=3$. It is observed that the relative deviation for the adaptive CFL is well controlled within bounds as expected. We compare schemes with $P^2$ and $P^3$ LDG solvers and observe that the scheme with $P^3$ LDG performs better in controlling relative deviation of upstream cell areas. On its subplot (b), we observe that quadratic-curved approximations to sides of upstream cells are crucial in preserving upstream cell areas. The scheme performs much better than the counterpart without the QC approximation. In fact, for this example, perhaps because the numerical mesh does not fully resolve the solution structures (thus some numerical oscillations appear), the $P^1$ SLDG scheme is performing better than the $P^2$ SLDG scheme in preserving upstream cell areas; while the $P^2$ SLDG-QC scheme performs the best.
 \item {\em Time histories of relative deviation of energy and enstrophy are plotted in Figure~\ref{figure:KH_norm}.} The performance of the SLDG schemes in preserving the invariants is comparable. We remark that, the proposed conservative SLDG schemes are able to better preserve energy than the non-conservative SLWENO scheme \cite{xiong}. We also note that, in many situations, higher order schemes preserve better these invariants; yet there are some exceptions which are subject to further investigation.
\end{enumerate}

\begin{figure}[h!]
\centering                              
\includegraphics[width=50mm, height=50mm]{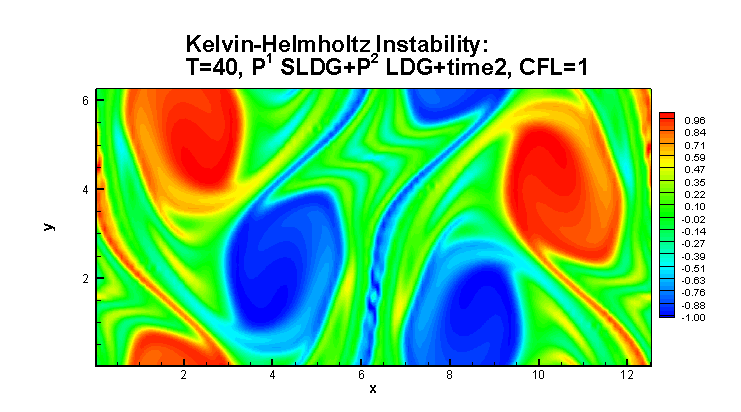}
\hspace{-2mm}
\includegraphics[width=50mm, height=50mm]{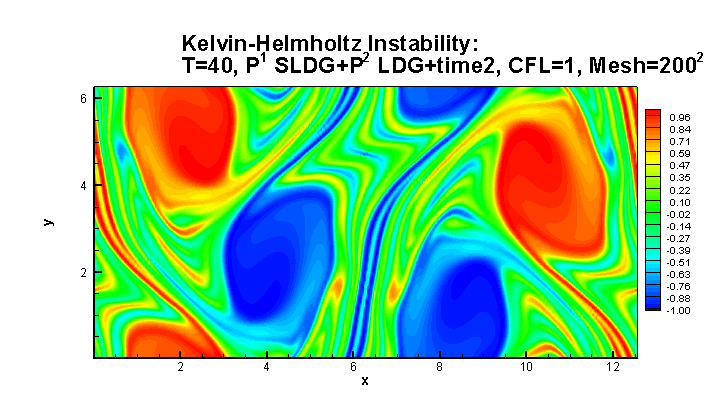}
\includegraphics[width=50mm, height=50mm]{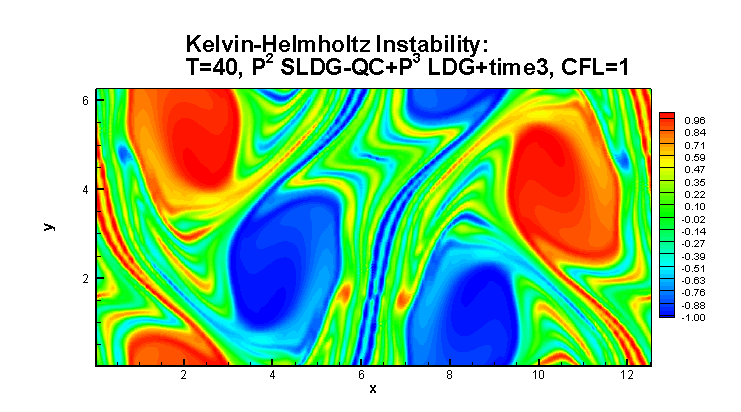}
\caption{Contour plots of the numerical solutions for the Kelvin-Helmholtz instability at $T=40$.
The mesh of $100\times100$ is used, unless otherwise specified.
}
\label{figure:KH_contour1}
\end{figure}

\begin{figure}[h!]
\includegraphics[width=79mm, height=50mm]{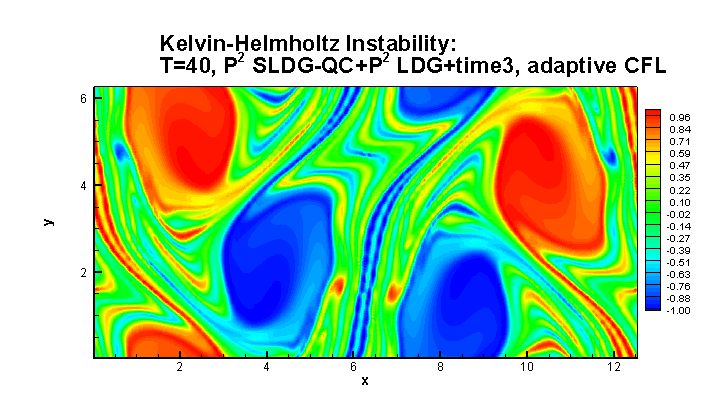}
\includegraphics[width=79mm, height=50mm]{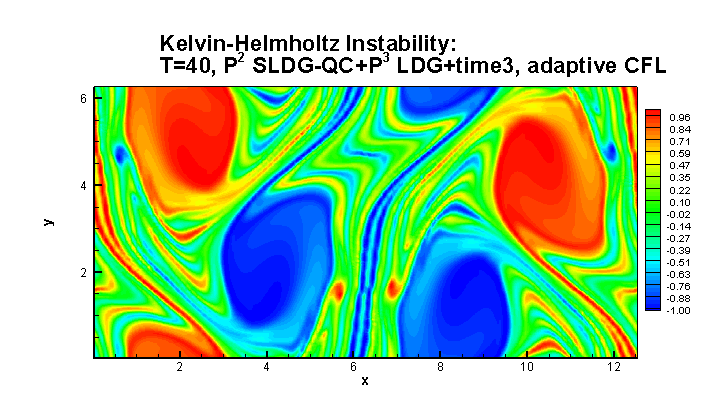}
\caption{Contour plots for the Kelvin-Helmholtz instability at $T=40$.
Schemes with adaptive CFL. $P^2$ LDG (left) and $P^3$ LDG (right) are used.
The mesh is $100\times100$.
}
\label{figure:KH_contour}
\end{figure}

\begin{figure}[h!]
\centering                              
\includegraphics[width=79mm]{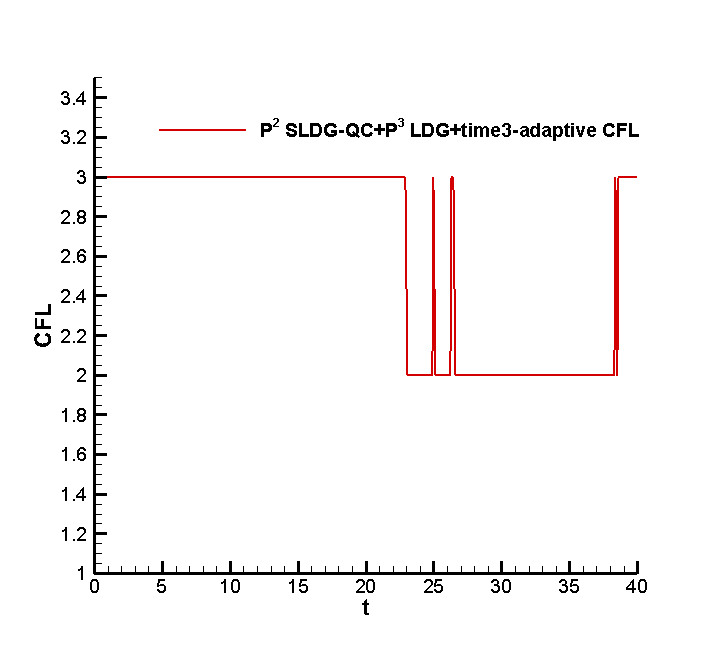}
\caption{The $CFL$ history of $P^2$ SLDG-QC+$P^3$ LDG+time3 with adaptive $CFL$ for the Kelvin-Helmholtz instability.
The mesh is $100\times100$.
}
\label{figure:KH_CFL}
\end{figure}

\begin{figure}[h!]
\centering
\subfigure[$P^2$ SLDG-QC versus $P^1$ SLDG]{
\includegraphics[width=75mm]{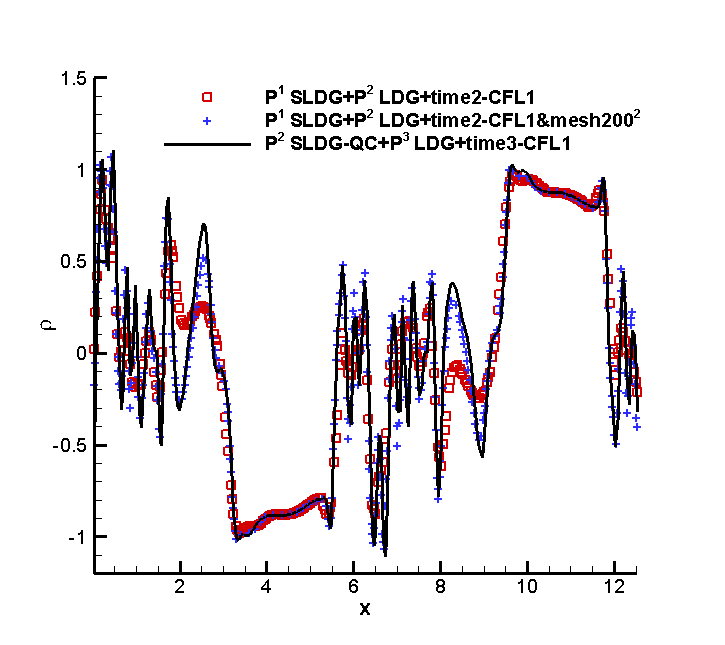}
}                         
\subfigure[$P^2$ SLDG-QC+$P^{3}$ LDG, adaptive $CFL$ versus non-adaptive $CFL$]{
\includegraphics[width=75mm]{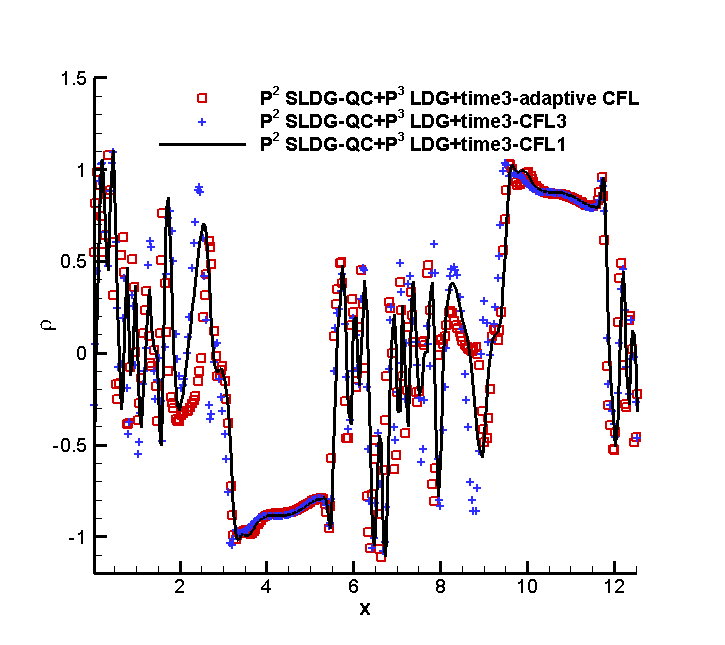}
}
\caption{1D cuts of the solutions at $y=\pi$ for the Kelvin-Helmholtz instability at $T=40$. The mesh of $100\times100$ is used unless otherwise specified. }
\label{figure:KH_1d}
\end{figure}

\begin{figure}[h!]
\centering
\subfigure[$P^2$ SLDG-QC]{
\includegraphics[width=75mm]{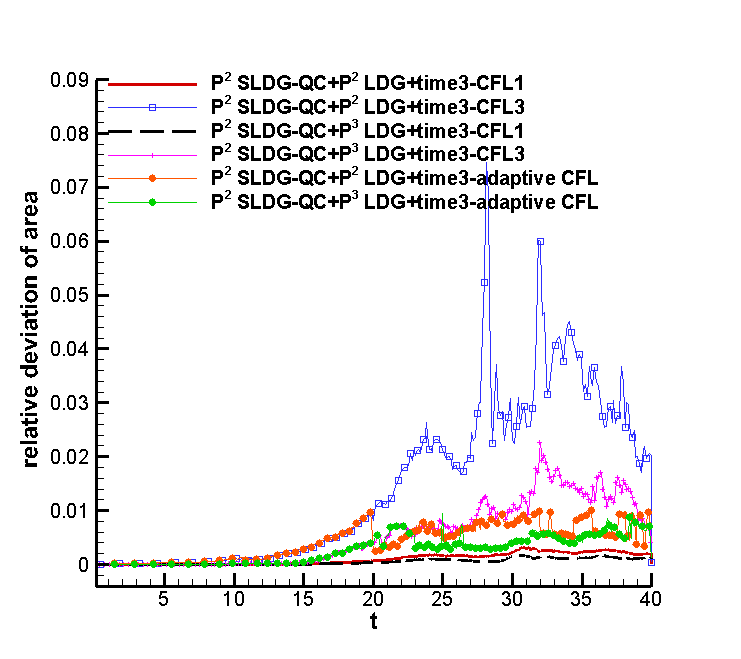}
}                         
\subfigure[$P^k$ SLDG-(QC) with $CFL=1$]{
\includegraphics[width=75mm]{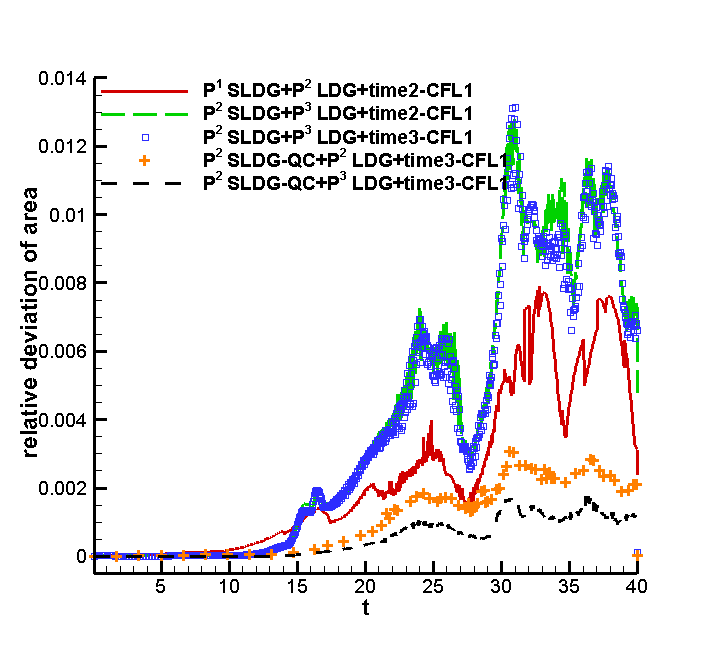}
}
\caption{Time evolution of the relative deviation of area for the proposed SLDG schemes for the Kelvin-Helmholtz instability. The mesh of $100\times100$ is used. }
\label{figure:KH_area}
\end{figure}

\begin{figure}[h!]
\centering
\subfigure[$P^2$ SLDG-QC]{
\includegraphics[width=75mm]{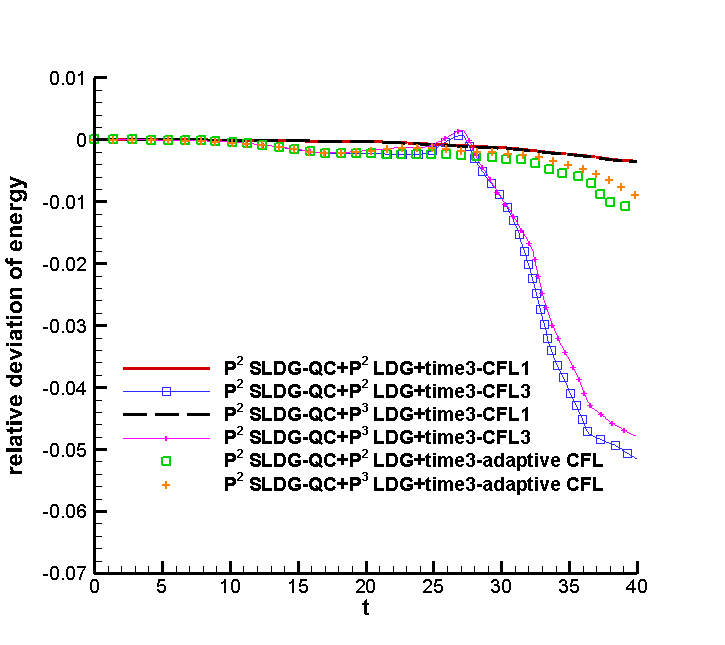}
}                         
\subfigure[$P^k$ SLDG-(QC) with $CFL=1$]{
\includegraphics[width=75mm]{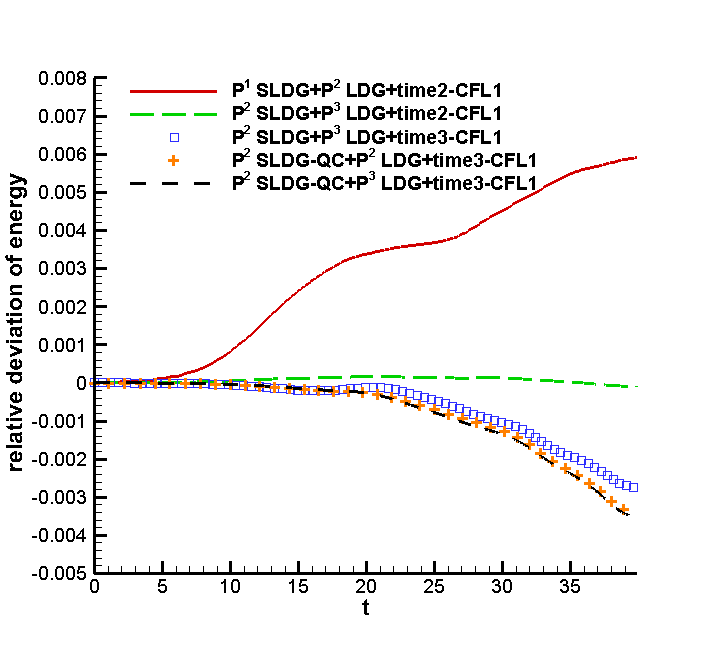}
}
\subfigure[$P^2$ SLDG-QC]{
\includegraphics[width=75mm]{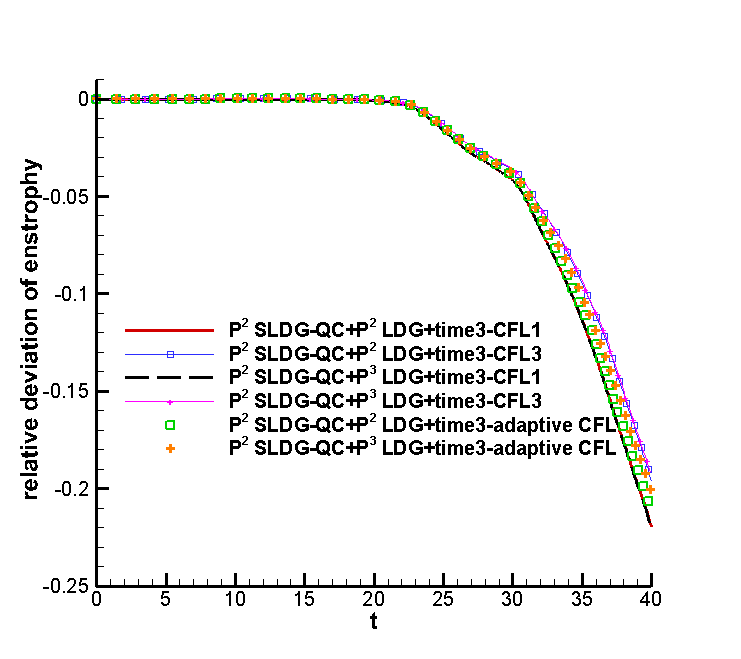}
}                         
\subfigure[$P^k$ SLDG-(QC) with $CFL=1$]{
\includegraphics[width=75mm]{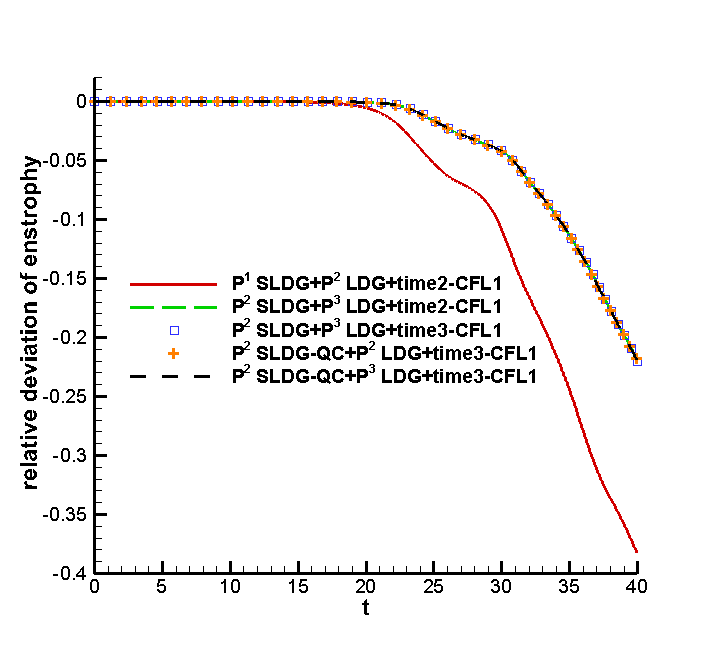}
}
\caption{Time evolution of the relative deviation of energy and enstrophy for the proposed SLDG methods for the Kelvin-Helmholtz instability. The mesh of $100\times100$ is used. }
\label{figure:KH_norm}
\end{figure}

\end{exa}

\begin{exa}
(The vortex patch problem) We solve the model problem \eqref{Euler} in the domain $[0,2\pi]\times[0,2\pi]$ with the initial condition
\begin{equation}
\omega(x,y,0) =
\begin{cases}
-1, & \text{if} \ (x,y) \in \left[\frac{\pi}{2},\frac{3\pi}{2}\right]\times\left[\frac{\pi}{4},\frac{3\pi}{4}\right], \\
1,  & \text{if} \ (x,y) \in \left[\frac{\pi}{2},\frac{3\pi}{2}\right]\times\left[\frac{5\pi}{4},\frac{7\pi}{4}\right],\\
0,  & \text{otherwise,}
\end{cases}
\end{equation}
and periodic boundary conditions.

When the proposed SLDG method with a large $CFL$ number, i.e. $CFL=3$, is applied, numerical oscillations are present. For example, in the left panel of Figure \ref{figure:vortex_WL}, we observe that the numerical solution computed by $P^2$ SLDG-QC+$P^3$ LDG+time3-CFL3 without the WENO limiter exhibits unphysical oscillatory behavior. When the WENO limiter is applied, numerical oscillations disappear, see the right panel of Figure \ref{figure:vortex_WL}.

When the SLDG scheme with adaptive $CFL$s is used, the initial $CFL=3$  is automatically reduced to $CFL=2$ at the beginning of the simulation due to the adaptive mechanism, see the right panel in Figure~\ref{figure:vortex_adaptive}. With adaptive $CFL$, it is observed that the scheme performs well in capturing solution structures without producing oscillations, even though no limiter is used. Extra robustness is observed from the adaptive time-stepping strategy. One intuitive explanation of the effect of ``controlling oscillations" by the adaptive CFL strategy is the following: without adaptive $CFL$, if a relatively large $CFL$ is used, numerical approximations to the shapes of upstream cells (and their areas) may not be accurate enough. Recall that when the areas of upstream cells are preserved, the maximum principle can be preserved (at least in terms of cell averages), oscillations can be avoided. When upstream cells are extremely distorted due to the large $CFL$ used, and then relative deviation of areas is likely to become larger. Consequently, cell averages of the solution may go out of bounds and become oscillatory. If no remedy is used, numerical approximations to upstream cells could become more distorted, leading to more pronounced unphysical oscillations.

As has been done before, we track the time history of the $L^\infty$ norm of relative deviation of areas of upstream cells, energy and enstrophy of various SLDG schemes in Figure \ref{figure:vortex_area} and Figure \ref{figure:vortex_norm}. The observation is similar to the previous example.

\begin{figure}[h!]
\centering                              
\includegraphics[width=75mm]{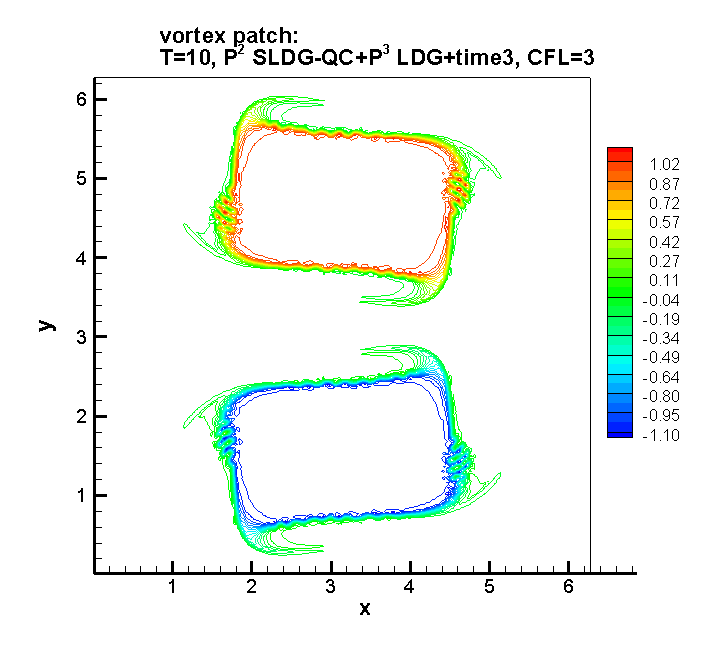}
\includegraphics[width=75mm]{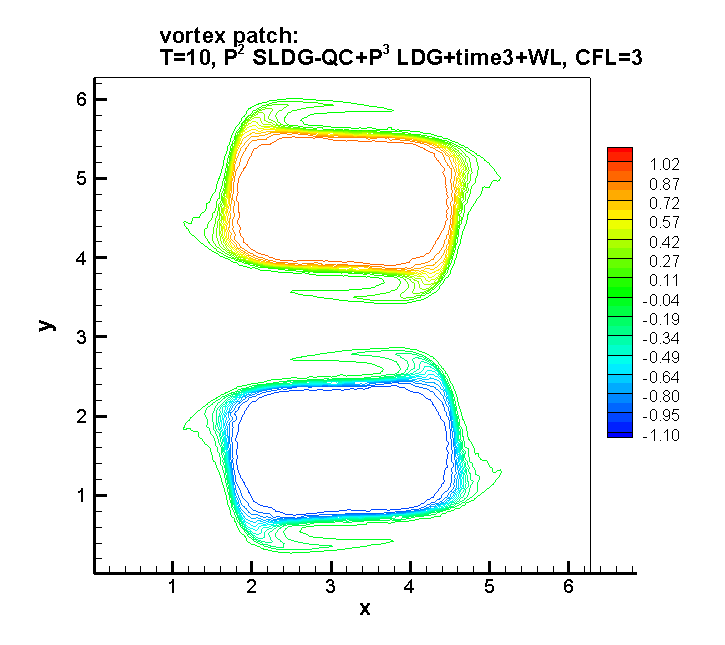}
\caption{Contour plots of the numerical solutions for $P^2$ SLDG-QC+$P^3$ LDG+time3 with (right) or without (left) WENO limiter for the vortex patch test.
The mesh of $100\times100$ is used.
30 equally spaced contours from $-1.1$ to $1.1$.
}
\label{figure:vortex_WL}
\end{figure}

\begin{figure}[h!]
\centering                              
\includegraphics[width=79mm]{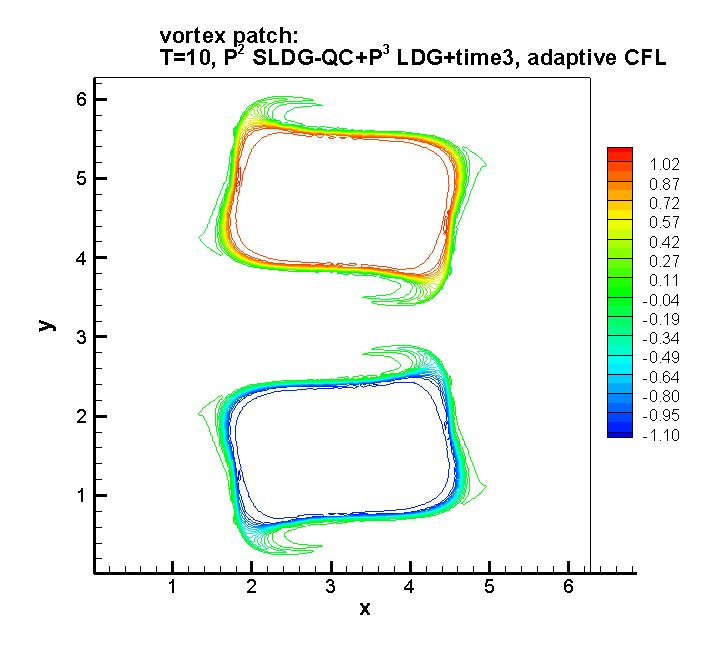}
\includegraphics[width=79mm]{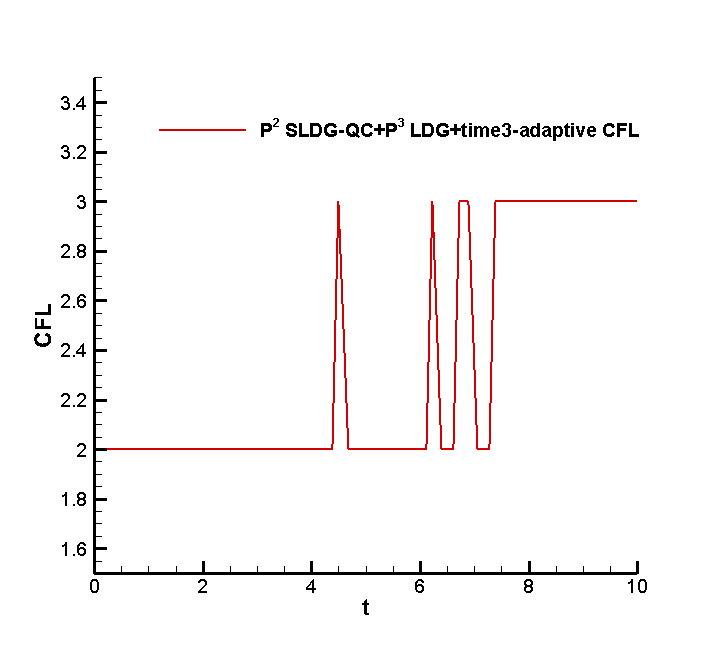}
\caption{
Left: The contour plot of the numerical solution for the vortex patch test solved by $P^2$ SLDG-QC+$P^3$ LDG+time3 with adaptive $CFL$;
right: time history of CFL. The mesh of $100\times100$ and the initial $CFL$ is $3$.
30 equally spaced contours from $-1.1$ to $1.1$.
}
\label{figure:vortex_adaptive}
\end{figure}

%
%

\begin{figure}[h!]
\centering
\subfigure[$P^2$ SLDG-QC]{
\includegraphics[width=75mm]{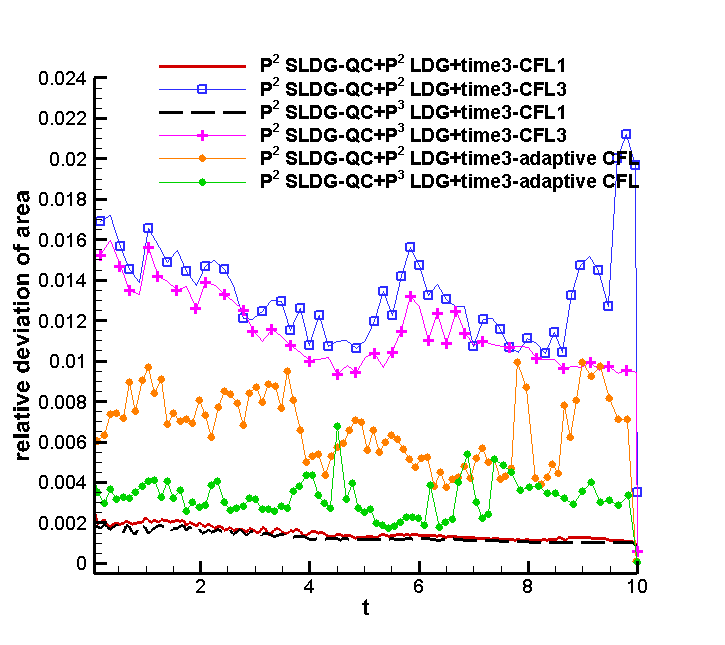}
}                         
\subfigure[$P^k$ SLDG-(QC) with $CFL=1$]{
\includegraphics[width=75mm]{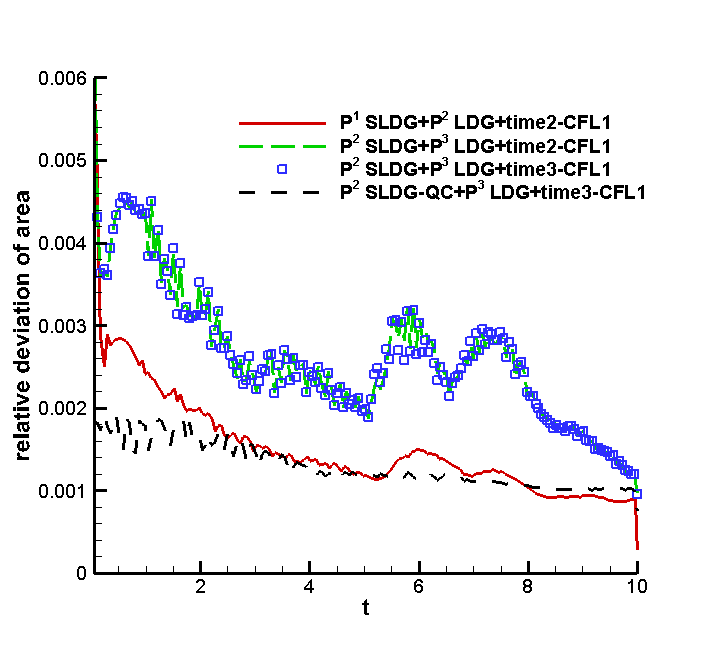}
}
\caption{Time evolution of the relative deviation of upstream areas for the proposed SLDG methods for the vortex patch test. The mesh of $100\times100$ is used. }
\label{figure:vortex_area}
\end{figure}

\begin{figure}[h!]
\centering
\subfigure[$P^2$ SLDG-QC]{
\includegraphics[width=75mm]{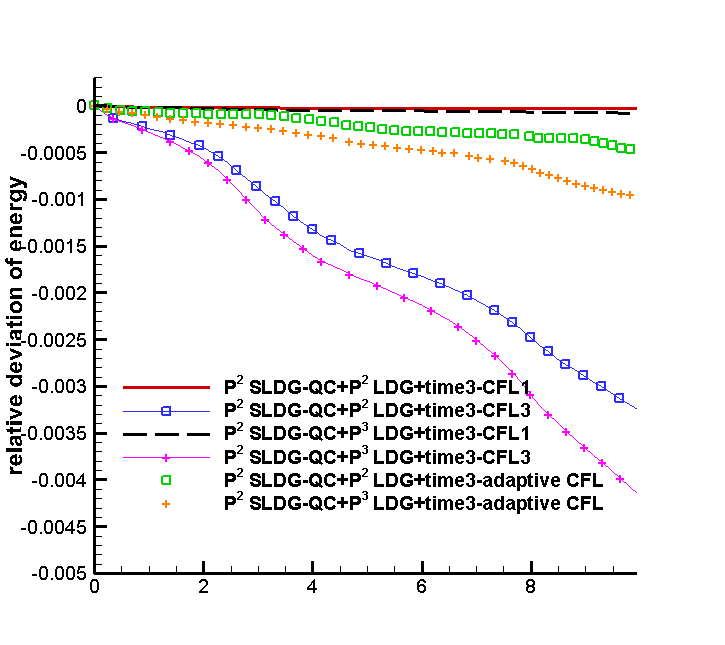}
}                         
\subfigure[$P^k$ SLDG-(QC) with $CFL=1$]{
\includegraphics[width=75mm]{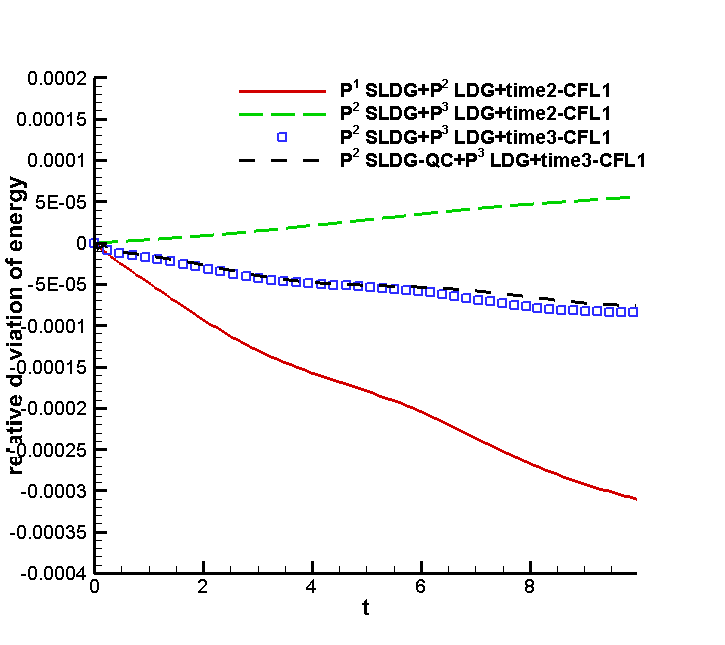}
}
\subfigure[$P^2$ SLDG-QC]{
\includegraphics[width=75mm]{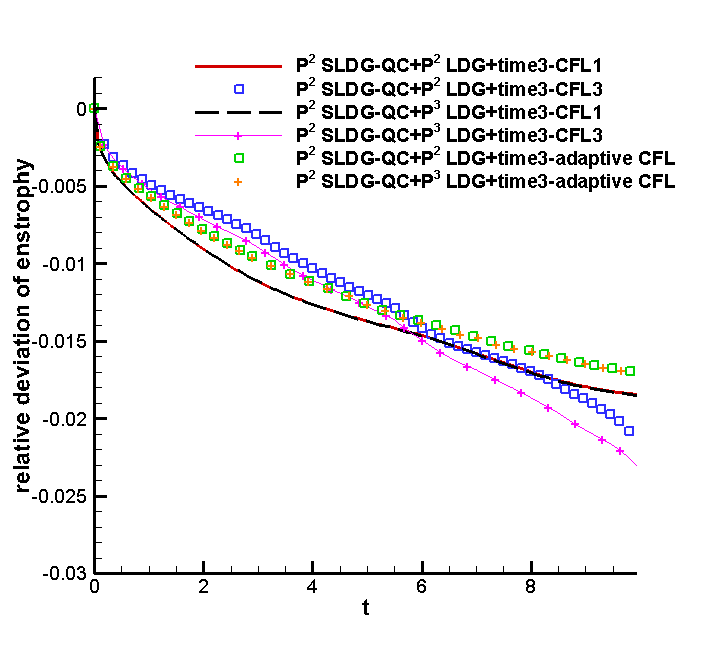}
}                         
\subfigure[$P^k$ SLDG-(QC) with $CFL=1$]{
\includegraphics[width=75mm]{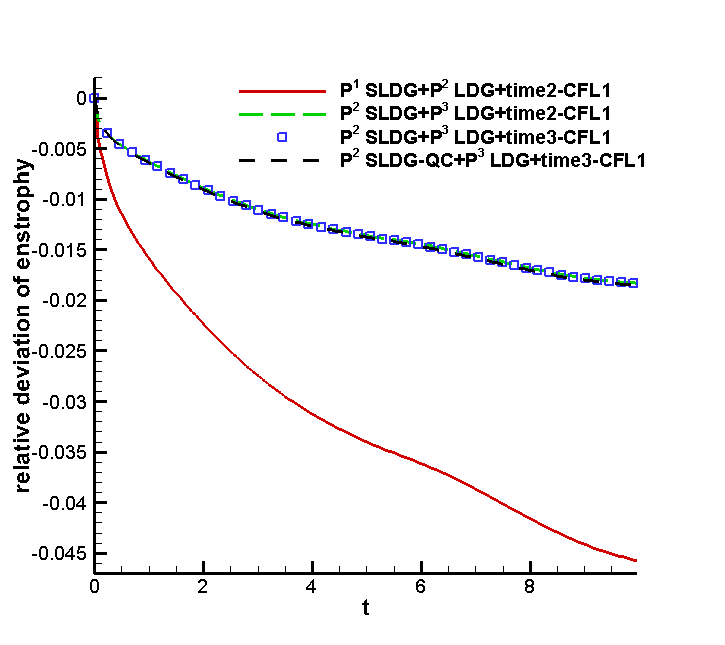}
}
\caption{Time evolution of the relative deviation of energy and enstrophy for the proposed SLDG methods for the vortex patch test. The mesh of $100\times100$ is used. }
\label{figure:vortex_norm}
\end{figure}


\end{exa}

\begin{exa}
(The shear flow problem)
For this double shear layer problem \cite{bell1989second,zhangshu2010}, we solve the model problem \eqref{Euler}
 in the domain $[0,2\pi]\times[0,2\pi]$, with periodic boundary conditions and
the initial condition given by
\begin{equation}
\omega(x,y,0)
=
\begin{cases}
\delta \cos(x) - \frac{1}{\rho} sech^2\left( \frac{y-\pi/2}{\rho} \right), & \text{if} \ y\leq \pi, \\
\delta \cos(x) + \frac{1}{\rho} sech^2\left( \frac{3\pi/2-y}{\rho} \right), & \text{if} \ y>\pi,
\end{cases}
\end{equation}
where $\delta =0.05$ and $\rho = \pi/15$.

As time evolves, the solution quickly develops into roll-ups with smaller and smaller spatial scales. On any fixed grid, the full resolution will be lost eventually. This problem has been tested by the high order Eulerian finite difference ENO/WENO method in \cite{weinan1994numerical,guo2015maximum}, the high order 	SLWENO scheme in \cite{qiu_shu_sl}, the DG method in \cite{liu2000high,zhangshu2010,zhu2017h} and the spectral element method in \cite{fischer2001filter,xu2006stabilization}.
We solve this problem up to $T=8$ by using the SLDG method with the mesh of $100\times100$ elements.
Figure \ref{figure:shear1} presents the solution for the shear flow test at $T=8$ solved by $P^2$ SLDG+$P^3$ LDG+time3 with $CFL=1$ without the WENO limiter (left) and with the WENO limiter (right). Numerical oscillations are observed (upper mid and lower right regions of the plot) for the scheme without the WENO limiter. Once the robust WENO limiter is applied, numerical oscillations disappear.
In Figure \ref{figure:shear_adaptive}, we show the contour plot of the solution computed by $P^2$ SLDG-QC+$P^3$ LDG+time3 with adaptive $CFL$ (left) as well as the $CFL$ history over time (right). No limiter is used, yet no oscillation is observed. Such results suggest that the adaptive $CFL$ time-stepping method improves the robustness of the SLDG schemes. The conclusion we draw in this example is similar to that in the vortex patch example.

We further compare the $L^\infty$ norm of relative deviation of upstream cell areas, as well as energy and enstrophy in Figure \ref{figure:shear_area} and Figure \ref{figure:shear_norm}, respectively. We compare the performance of the scheme with adaptive $CFL$ without the WENO limiter, with that from the SLDG schemes with a fixed $CFL$ and with or without the WENO limiter. It is observed that the scheme without WENO limiter performs better than that with the WENO limiter (even with smaller CFL number) in preserving the invariants. Such a phenomenon can be explained by the fact that, when the solution is under-resolved, the WENO limiter is often turned on. By using lower order polynomials in the approximation space, larger deviations may be induced.


\begin{figure}[h!]
\centering                              
\includegraphics[width=75mm]{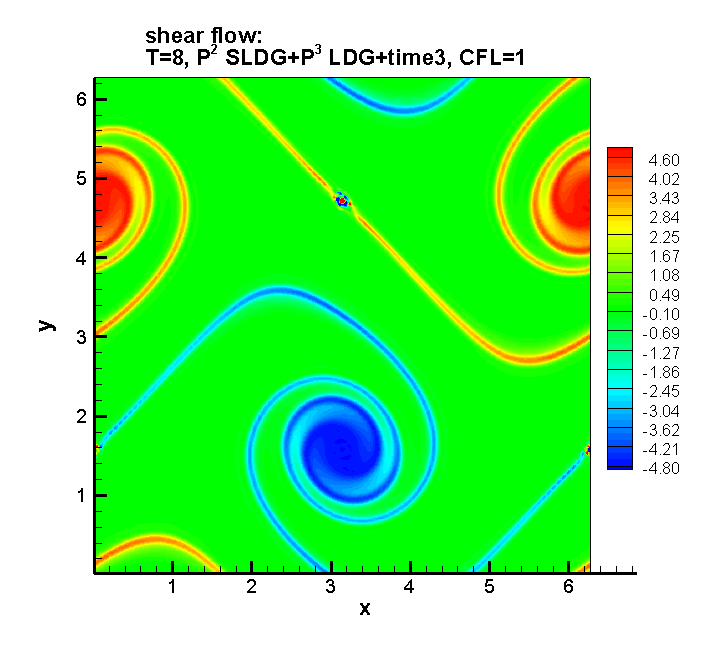}
\includegraphics[width=75mm]{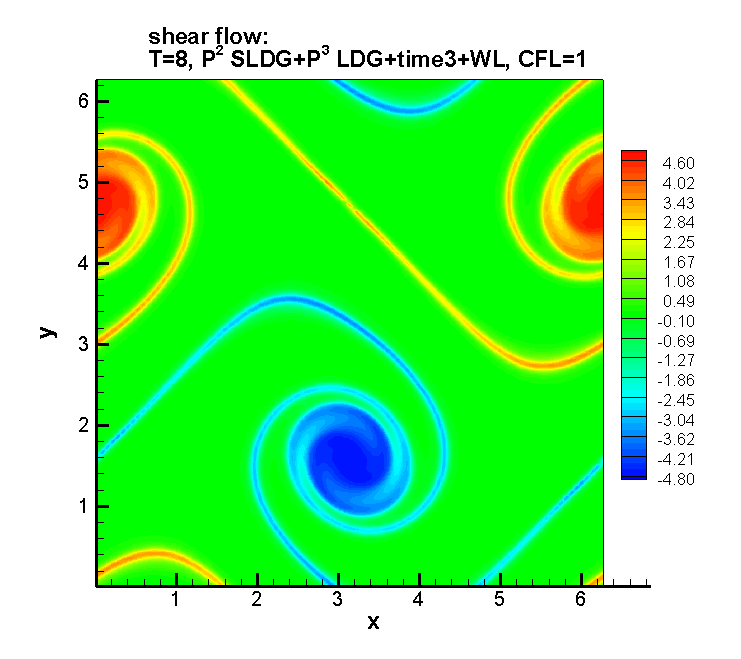}
\caption{Contour plots of the numerical solutions for the shear flow test at $T=8$ solved by the methods without WENO limiter (left) and with WENO limiter (right).
The mesh of $100\times100$ and $CFL=1$.
}
\label{figure:shear1}
\end{figure}

\begin{figure}[h!]
\centering                              
\includegraphics[width=79mm]{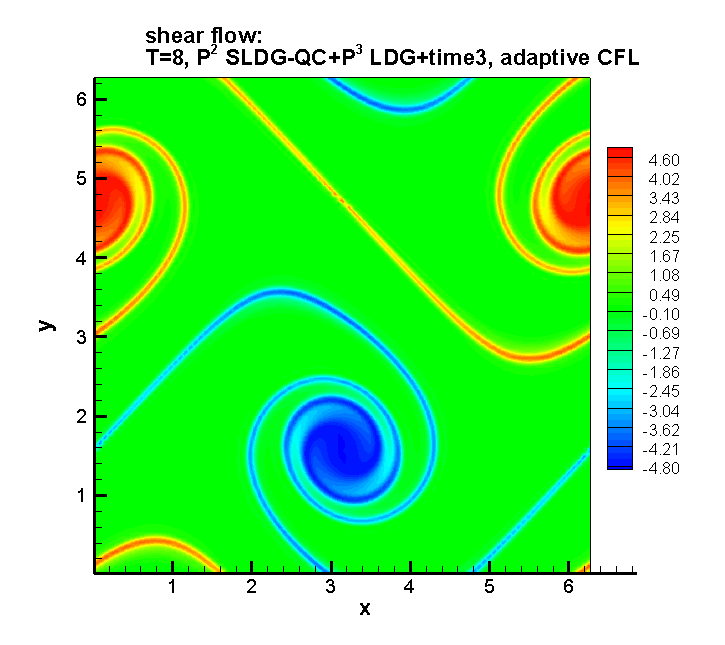}
\includegraphics[width=79mm]{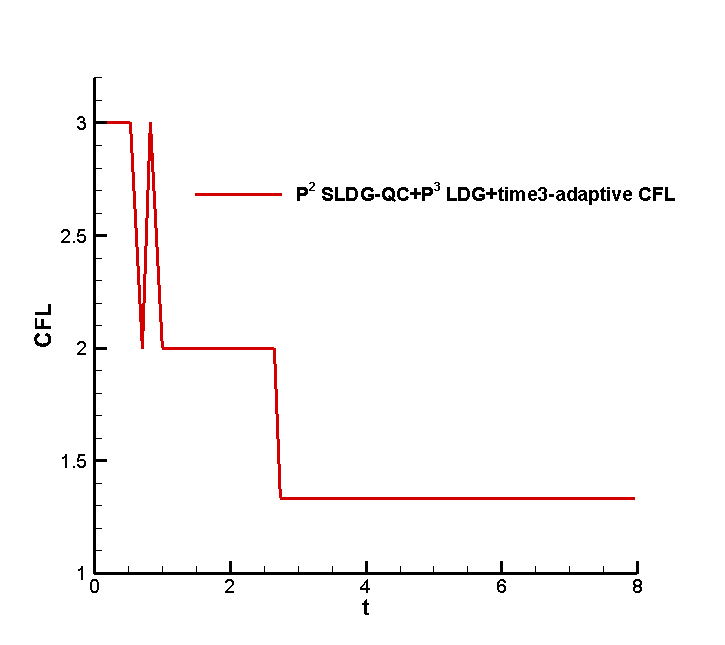}
\caption{Left: The contour plot of the numerical solution for the shear flow test solved by $P^2$ SLDG-QC+$P^3$ LDG+time3 with adaptive $CFL$;
right: time history of CFL. The mesh of $100\times100$ and the initial $CFL$ is $3$.
}
\label{figure:shear_adaptive}
\end{figure}


\begin{figure}[h!]
\centering
\subfigure[$P^2$ SLDG-QC]{
\includegraphics[width=75mm]{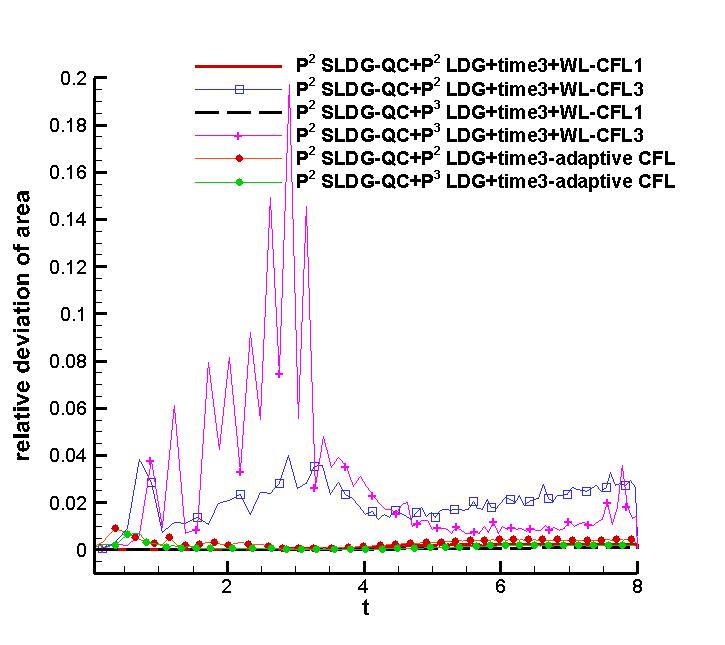}
}                         
\subfigure[$P^k$ SLDG-(QC) with $CFL=1$]{
\includegraphics[width=75mm]{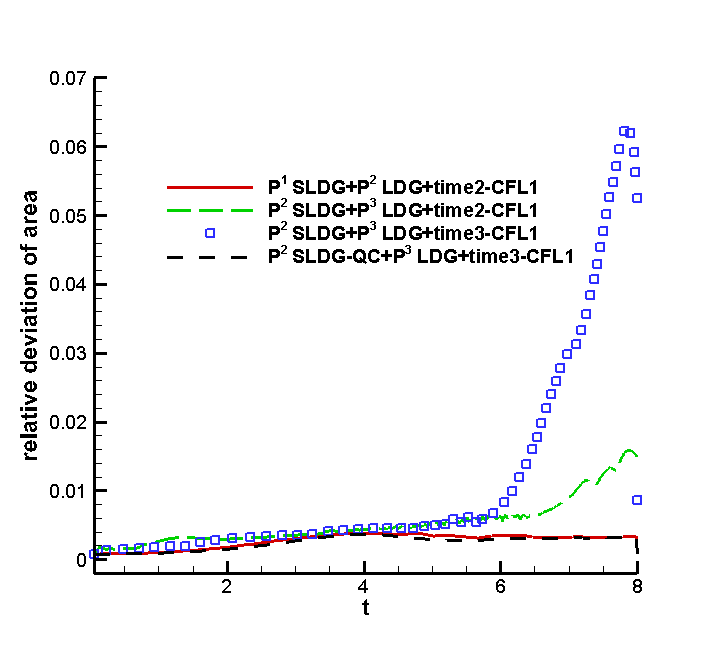}
}
\caption{Time evolution of the relative deviation of areas of upstream cells for the proposed SLDG scheme for the shear flow test. The mesh is $100\times100$. }
\label{figure:shear_area}
\end{figure}

\begin{figure}[h!]
\centering
\subfigure[$P^2$ SLDG-QC]{
\includegraphics[width=75mm]{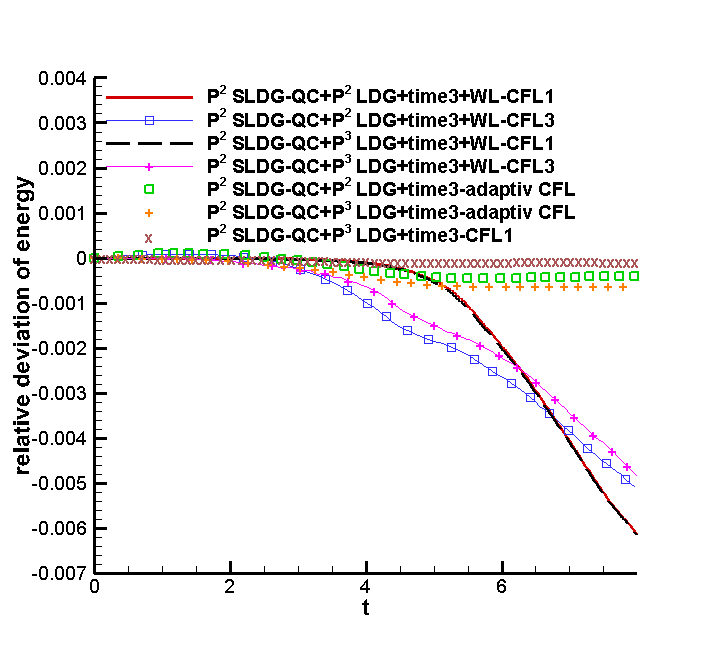}
}                         
\subfigure[$P^k$ SLDG-(QC) with $CFL=1$]{
\includegraphics[width=75mm]{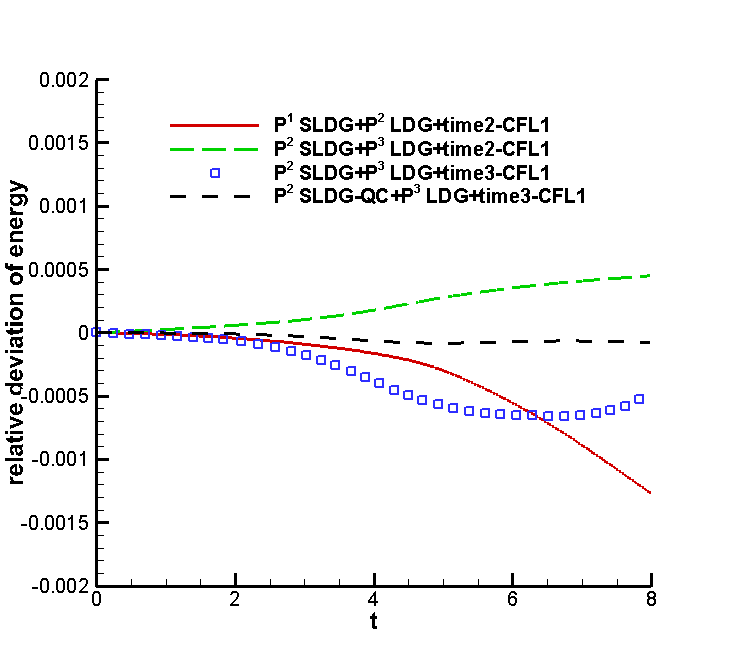}
}
\subfigure[$P^2$ SLDG-QC]{
\includegraphics[width=75mm]{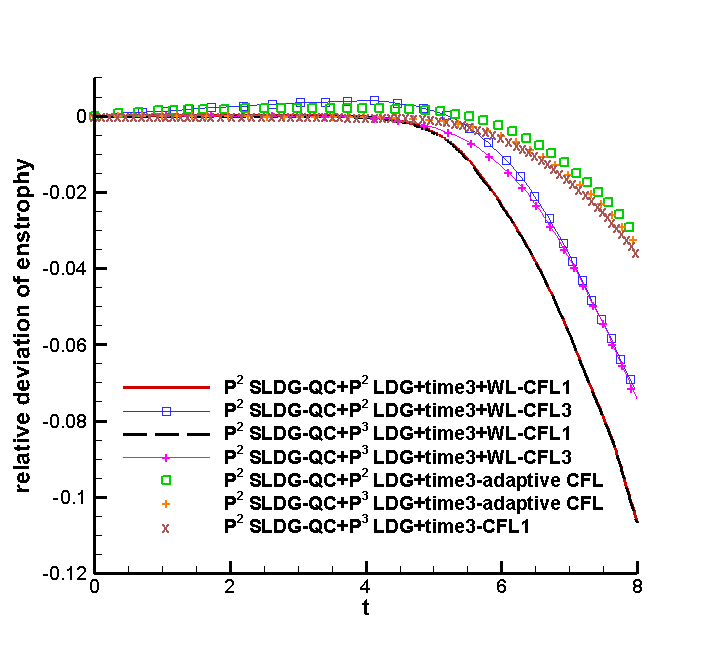}
}                         
\subfigure[$P^k$ SLDG-(QC) with $CFL=1$]{
\includegraphics[width=75mm]{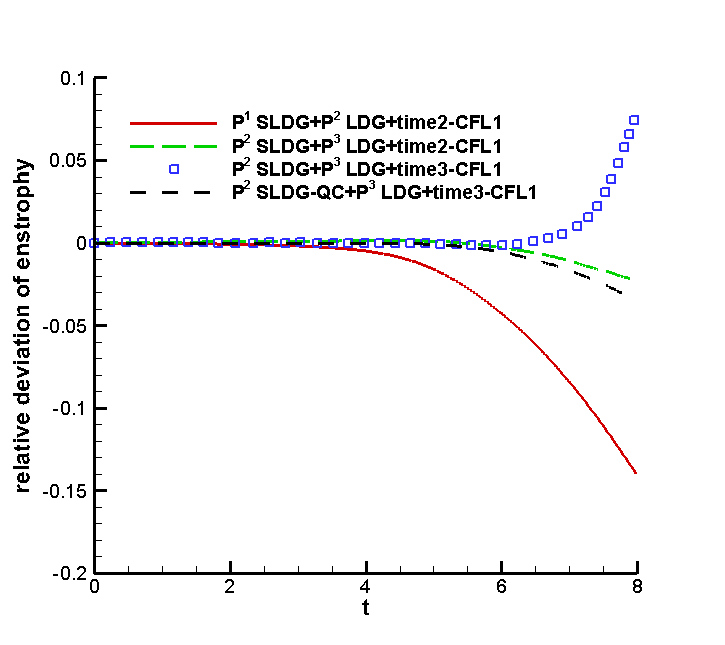}
}
\caption{Time evolution of the relative deviation of energy and enstrophy for the proposed SLDG schemes for the shear flow test. The mesh is $100\times100$. }
\label{figure:shear_norm}
\end{figure}
\end{exa}

\section{Conclusion}
\label{sec6}
\setcounter{equation}{0}
\setcounter{figure}{0}
\setcounter{table}{0}

In this paper, we proposed a high order conservative semi-Lagrangian discontinuous Galerkin (SLDG) method for solving two-dimensional incompressible Euler equations and the guiding center Vlasov model without operator splitting. The three key ingredients include a high order conservative SLDG transport scheme as the backbone of the algorithm, a high order characteristics tracing technique, and an adaptive time-stepping strategy to further enhance the robustness and effectiveness of the scheme. As the main advantage, the scheme is able to take large time step evolution, and at the same time be high order accurate in both space and time and mass conservative.  The performance of the scheme in terms of order accuracy in space and time, CPU cost as well as the ability to preserve important physical invariants was benchmarked though extensive numerical experiments. We only consider periodic boundary condition in this paper. The extension to general boundary conditions is subject to our future work.


\bibliographystyle{abbrv}
\bibliography{refer17}

\end{document}